\documentclass[12pt,twoside,a4paper,reqno]{amsart}
\usepackage[a4paper,width=160mm,top=25mm,bottom=30mm]{geometry}
\usepackage{amsmath,amsthm,amssymb}
\usepackage{hyperref}
\usepackage{amsfonts}
\usepackage{enumerate}
\usepackage[normalem]{ulem}
\usepackage{color}

\newtheorem{theorem}{Theorem}

\newtheorem{conjecture}[theorem]{Conjecture}
\newtheorem{corollary}[theorem]{Corollary}

\newtheorem{definition}[theorem]{Definition}

\newtheorem{lemma}[theorem]{Lemma}

\newtheorem{proposition}[theorem]{Proposition}
\newtheorem{remark}[theorem]{Remark}

\input{tcilatex}

\begin{document}
\title{Gaussian fluctuations for the classical XY model}
\author{Charles M. Newman, Wei Wu}
\address[Charles Newman]{Courant Institute of Mathematical Sciences, New
York University, 251 Mercer st, New York, NY 10012, USA \\
\& NYU-ECNU Institute of Mathematical Sciences at NYU Shanghai, 3663
Zhongshan Road North, Shanghai 200062, China.}
\address[Wei Wu]{Courant Institute of Mathematical Sciences, New York
University, 251 Mercer st, New York, NY 10012, USA \\
\& NYU-ECNU Institute of Mathematical Sciences at NYU Shanghai, 3663
Zhongshan Road North, Shanghai 200062, China.}

\begin{abstract}
We study the classical XY model in bounded domains of $\mathbb{Z}^{d}$ with
Dirichlet boundary conditions. We prove that when the temperature goes to
zero faster than a certain rate as the lattice spacing goes to zero, the
fluctuation field converges to a standard Gaussian white noise. This and
related results also apply to a large class of gradient field models.
\end{abstract}

\maketitle

\section{Introduction}

In this paper we study the classical XY model (also known as the classical
rotator model) in dimension $d\geq 2$. Assign to each vertex $i\in \mathbb{Z}%
^{d}$ a spin $s\left( i\right) \in S^{1}$ with the corresponding angle $%
\theta \left( i\right) \in \lbrack -\pi ,\pi )$. The XY\ model is defined
formally as a Gibbs measure with Hamiltonian given by%
\begin{equation*}
H=-\sum_{i\sim j}s\left( i\right) \cdot s\left( i\right) =-\sum_{i\sim
j}\cos \left( \theta \left( i\right) -\theta \left( j\right) \right) ,
\end{equation*}%
where the sum is over all nearest neighbor pairs of vertices. It is an
example of the more general $O\left( n\right) $ model, where each spin takes
a value in the sphere $S^{n-1}$.

A simple heuristic discussion of the low temperature behavior of the XY
model is as follows. As the temperature goes to zero, the spins tend to
align with each other so as to minimize the Hamiltonian. Since $\cos \left(
\delta \theta \right) \approx 1-\left( \delta \theta \right) ^{2}/2$ for
small $\delta \theta $, it is expected that at low temperature, the XY Gibbs
measure on large scales behaves like a Gaussian Free Field (GFF). This idea
originated in \cite{Dy} (see also \cite{MW}) and was referred to as the
Gaussian spin-wave approximation. By further making connections between
rotational symmetry of the XY model and the recurrence/transience property
of simple random walks, it was proved in \cite{MW} that for $d\leq 2$, there
is no spontaneous magnetization at any strictly positive temperature. A
related argument was applied in \cite{FSS} to show that for $d\geq 3$, with
appropriate boundary conditions, the spin configuration has a preferred
direction at low temperature --- i.e., there is spontaneous magnetization,
but has no spontaneous magnetization at high temperature (and thus there is
a phase transition).

For $d=2$, the Gaussian spin-wave approximation is expected to be valid for
all temperatures below the critical temperature corresponding to the
so-called Kosterlitz-Thouless transition (see \cite{KT}). This suggests, for
example, that the scaling limits, as the lattice spacing tends to zero, of
the XY and GFF models on the lattice, should be closely related to each
other. But this has not been proved for any fixed positive temperature. The
best known results are only for the two point correlation functions: a
polynomial upper bound for the spin-spin correlation function was proved in 
\cite{MS}, and a different polynomial lower bound was established in \cite%
{FSp}. A more modest question, which is to show that at some fixed low
temperature, the spin-spin correlation has an asymptotic power law decay,
still remains open. For $d\geq 3$, it was proved that at low temperature the
spin configuration has long range order by using the infrared bound \cite%
{FSS} and also proved that the two point (transverse) correlation function
has a power law decay \cite{BFLLS}. Yet it is still an open problem to match
the powers in the upper and lower bound for the longitudinal correlation
function (the cosine-cosine correlation).

In this paper we make some modest progress about the Gaussian spin-wave
approximation. We study the XY model in finite domains in $\varepsilon 
\mathbb{Z}^{d}$, for $d\geq 2$, with Dirichlet boundary conditions, where
the inverse temperature $\beta $ depends on $\varepsilon $, and is required
to grow at least proportional to $|\log \varepsilon |$ as $\varepsilon $
tends to zero (see Section \ref{model} for precise conditions). With this
assumption, the spin configuration has a preferred direction even for $d=2$,
and our results focus on the fluctuation field. We prove that the rescaled
(properly defined) gradient field of the XY model converges weakly to a
Gaussian white noise. This suggests that the spin field should converge to
some version of a compactified GFF (i.e., a GFF modulo $2\pi $).

The significance of our requirement on $\beta \left( \varepsilon \right) $
growing fast enough as $\varepsilon \rightarrow 0$ is that when $\beta
\left( \varepsilon \right) $ is proportional to $|\log \varepsilon |$,
applying a result of \cite{BF} (see our Theorem \ref{xypierels}) shows that
on an event with high probability, there are no vortices (see Section \ref%
{model} for the definition) in a domain of size $\varepsilon ^{-1}$. This
allows one to construct a coupling between the XY model and a corresponding
convex gradient field model, such that their gradient variables agree on an
event with high probability. It then relates the scaling limit of the XY
model to the corresponding limit for a gradient field with a convex
potential, where the situation is much better understood because of the
Helffer-Sj{\"{o}}strand formula (\cite{HS}, see also our Section \ref{HS})
that represents the covariance of the local functionals of the field in
terms of a random walk in random conductances. We mention \cite{NS}, \cite%
{GOS} and \cite{M} as an incomplete list of literature regarding central
limit theorem type results for such convex gradient fields. We adapt the
arguments in these papers to the current setting where the potential depends
on the lattice spacing, in order to prove our fluctuation results. We
further show that the limiting quadratic form has intensity one --- i.e.,
the white noise scaling limit corresponds to the gradient of a standard GFF.

The rest of this paper is organized as follows. In Section \ref{model} we
set up the model and state the main results (Theorems \ref{dirichletXY} and %
\ref{gf}). Section \ref{couple} starts with a contour estimate for the XY
model (Theorem \ref{xypierels}) from \cite{BF}, and a corresponding result
for the gradient field (Proposition \ref{convexpierels}); it then uses these
estimates to construct a coupling between the two measures. The rest of the
paper is then devoted to central limit theorem (CLT) type results for the
corresponding gradient field model. We first prove a CLT type result in all
of $\varepsilon \mathbb{Z}^{d}$ (in Section \ref{Zd}) and then in bounded
domains (in Section \ref{bdd}). Combined with the coupling established in
Section \ref{couple} we complete the proof of our main results. Finally,
Section \ref{open} discusses two open questions and related conjectures.

\section{Main Results\label{model}}

Fix $d\geq 2$. Let $D\subset \mathbb{R}^{d}$ be a bounded domain with smooth
boundary, $D^{\varepsilon }=D\cap \varepsilon \mathbb{Z}^{d}$. We also
denote by $(D^{\varepsilon })^{\ast }$ the set of all nearest neighbor
directed edges $\left( x,y\right) $ with $x$ and $y$ in $D^{\varepsilon }$.
We will have a further assumption on the $(D^{\varepsilon })^{\ast }$ graph
and hence on $D$ below. Given $x\in D^{\varepsilon }$, we associate an
angular variable $\theta ^{\varepsilon }\left( x\right) \in \lbrack -\pi
,\pi )$, and a spin variable $s^{\varepsilon }\left( x\right) =\left( \cos
\theta ^{\varepsilon }\left( x\right) ,\sin \theta ^{\varepsilon }\left(
x\right) \right) $. Given $\beta >0$, the classical XY model on $%
D^{\varepsilon }$ is defined by the following Gibbs measure $\mu _{\beta
}^{\varepsilon }$ on $[-\pi ,\pi )^{D^{\varepsilon }}$ (or $\left(
S^{1}\right) ^{D^{\varepsilon }}$):%
\begin{equation}
d\mu _{\beta }^{\varepsilon }=Z_{\beta }^{-1}\exp \left[ \beta \sum_{\left(
i,j\right) \in \left( D^{\varepsilon }\right) ^{\ast }}\cos \left( \theta
^{\varepsilon }\left( i\right) -\theta ^{\varepsilon }\left( j\right)
\right) \right] \prod_{i\in D^{\varepsilon }\backslash \partial
D^{\varepsilon }}d\theta ^{\varepsilon }\left( i\right) \prod_{i\in \partial
D^{\varepsilon }}\delta _{0}\left( d\theta ^{\varepsilon }\left( i\right)
\right) ,  \label{xygibbs}
\end{equation}%
where $Z_{\beta }$ is the normalizing constant. For any directed edge $%
b=\left( i,j\right) \in (D^{\varepsilon })^{\ast }$, we define the variable 
\begin{equation}
\eta ^{\varepsilon }\left( b\right) =\left\{ 
\begin{array}{cc}
\theta ^{\varepsilon }\left( j\right) -\theta ^{\varepsilon }\left( i\right)
+2\pi & \text{if }\theta ^{\varepsilon }\left( j\right) -\theta
^{\varepsilon }\left( i\right) \in \left( -2\pi ,-\pi \right) \\ 
\theta ^{\varepsilon }\left( j\right) -\theta ^{\varepsilon }\left( i\right)
& \text{if }\theta ^{\varepsilon }\left( j\right) -\theta ^{\varepsilon
}\left( i\right) \in \lbrack -\pi ,\pi ) \\ 
\theta ^{\varepsilon }\left( j\right) -\theta ^{\varepsilon }\left( i\right)
-2\pi & \text{if }\theta ^{\varepsilon }\left( j\right) -\theta
^{\varepsilon }\left( i\right) \in \lbrack \pi ,2\pi )%
\end{array}%
\right. .  \label{eta}
\end{equation}%
In other words, $\eta ^{\varepsilon }\left( b\right) $ is the smallest
increment (in absolute value) that deforms $\theta ^{\varepsilon }\left(
i\right) $ to $\theta ^{\varepsilon }\left( j\right) $ ($\func{mod}2\pi $).
For $l=1,...,d$ and $x\in D^{\varepsilon }$, we will sometimes denote $\eta
_{l}^{\varepsilon }\left( x\right) =\eta ^{\varepsilon }\left(
x,x+\varepsilon e_{l}\right) $, where $e_{l}$ is the $l^{th}$ standard basis
vector in $\mathbb{R}^{d}$ or $\mathbb{Z}^{d}$. The Hamiltonian\ of the XY
model (more precisely, $\beta $ times the Hamiltonian) can also be written
as $H\left( \eta \right) =-\beta \sum_{b\in \left( D^{\varepsilon }\right)
^{\ast }}\cos \left( \eta ^{\varepsilon }\left( b\right) \right) $.

We see that $\left\{ \eta ^{\varepsilon }\right\} $ satisfies the following
hardcore constraint:%
\begin{equation}
\eta ^{\varepsilon }\left( b_{1}\right) +\eta ^{\varepsilon }\left(
b_{2}\right) +\eta ^{\varepsilon }\left( b_{3}\right) +\eta ^{\varepsilon
}\left( b_{4}\right) =2k_{P}\pi ,  \label{kp}
\end{equation}%
whenever $\left( b_{1},b_{2},b_{3},b_{4}\right) $ encloses a plaquette $P$
in $D^{\varepsilon }$, where $k_{P}=0,\pm 1,\pm 2$. If $k_{P}\neq 0$, we
call $P$ a vortex of the spin configuration $\left\{ \theta ^{\varepsilon
}\right\} $. $\eta ^{\varepsilon }$ differs from the usual (discrete)
gradient of $\theta ^{\varepsilon }$, which is just $\theta ^{\varepsilon
}\left( j\right) -\theta ^{\varepsilon }\left( i\right) $ for $\left(
i,j\right) \in (D^{\varepsilon })^{\ast }$. When $\left\{ \theta
^{\varepsilon }\right\} $ is "vortex free", i.e., $\eta ^{\varepsilon
}\left( b_{1}\right) +\eta ^{\varepsilon }\left( b_{2}\right) +\eta
^{\varepsilon }\left( b_{3}\right) +\eta ^{\varepsilon }\left( b_{4}\right)
=0$ for every plaquette $P$, we want $\eta ^{\varepsilon }$ to be gradient
of some function $\Phi $, such that for any $x\in D^{\varepsilon }$, $\theta
^{\varepsilon }\left( x\right) =\Phi \left( x\right) \mod 2\pi $.

For the vortex free condition to imply the existence of such a $\Phi $, one
wishes to define $\Phi \left( x\right) =\sum_{b\in \mathcal{C}_{x^{\ast
}\rightarrow x}}\eta ^{\varepsilon }\left( b\right) ,$ where $\mathcal{C}%
_{x^{\ast }\rightarrow x}$ is any path in $(D^{\varepsilon })^{\ast }$
connecting $x^{\ast }$ to $x$ and conclude that $\eta ^{\varepsilon }=\nabla
\Phi $. For this to work it is necessary that $\Phi $ is well defined, i.e.,
that the vortex condition implies that the formula for $\Phi $ does not
depend on the choice of $\mathcal{C}_{x^{\ast }\rightarrow x}$. This will be
so if $(D^{\varepsilon })^{\ast }$ satisfies the following topological
requirement.

\begin{definition}
We say that $(D^{\varepsilon })^{\ast }$ is \textbf{simply connected} if for
any site self-avoiding loop $b$, there is a finite sequence of loops $%
b_{1}=b,$ $b_{2},...,b_{N}$ with $b_{j}$ and $b_{j-1}$ differing
(algebraically) by a single plaquette for each $j$, such that $b_{N}$ is a
trivial loop with no edge. We say that the domain $D$ in $\mathbb{R}^{d}$ is 
\textbf{discretely simply connected} if for each $\varepsilon >0$ such that $%
(D^{\varepsilon })^{\ast }$ is nonempty, $(D^{\varepsilon })^{\ast }$ is
simply connected.
\end{definition}

Henceforth we assume that $D$ is discretely simply connected. This, for
example, is the case for $D$ a rectangular parallelipiped along the
coordinate axis, $\left[ a_{1},b_{1}\right] \times ...\times \left[
a_{d},b_{d}\right] $ with $a_{j}<b_{j}$ for each $j$.

We study the zero temperature limit of the Gibbs measure (\ref{xygibbs}), by
letting $\beta =\beta \left( \varepsilon \right) $ depend on the lattice
size $\varepsilon $ with $\beta \left( \varepsilon \right) \rightarrow
\infty $ as $\varepsilon \rightarrow 0$. More precisely, we assume that, 
\begin{equation}
\beta \left( \varepsilon \right) +9d\log \varepsilon \rightarrow \infty 
\text{ as }\varepsilon \rightarrow 0.  \label{temp}
\end{equation}%
It is straightforward to see that with this choice of $\beta \left(
\varepsilon \right) $, the spins $\left\{ \theta ^{\varepsilon }\left(
x\right) \right\} _{x\in D^{\varepsilon }}$ (and the associated gradient
variables $\eta ^{\varepsilon }\left( b\right) $, $b\in \left(
D^{\varepsilon }\right) ^{\ast }$) concentrate around $0$ as $\varepsilon $
tends to zero (for a quantitative bound, see Theorem \ref{xypierels} below).
We focus on the fluctuation field of the $\left\{ \eta ^{\varepsilon
}\right\} $ variables. To do this, for $b\in \left( D^{\varepsilon }\right)
^{\ast }$, we define rescaled variables by setting $\tilde{\eta}%
^{\varepsilon }\left( b\right) =\sqrt{\beta \left( \varepsilon \right) }\eta
^{\varepsilon }\left( b\right) $. \bigskip Our main result is the following.

\begin{theorem}
\label{dirichletXY}Suppose that $\tilde{\eta}^{\varepsilon }=\sqrt{\beta
\left( \varepsilon \right) }\eta ^{\varepsilon }$ is the fluctuating field
associated with the Gibbs measure (\ref{xygibbs}), and suppose $\beta \left(
\varepsilon \right) $ satisfies (\ref{temp}). Define the linear functional 
\begin{equation}
\left\langle \tilde{\eta}^{\varepsilon },\varphi \right\rangle =\varepsilon
^{\left( d/2\right) -1}\sum_{b\in \left( D^{\varepsilon }\right) ^{\ast
}}\nabla \varphi \left( b\right) \tilde{\eta}^{\varepsilon }\left( b\right) 
\text{, \ for }\varphi \in \mathcal{C}_{0}^{\infty }\left( D\right) \text{,}
\label{xid}
\end{equation}%
where $\nabla \varphi \left( b\right) =\varphi \left( y\right) -\varphi
\left( x\right) $ for $b=\left( x,y\right) $. Then for any $t\in \mathbb{R}$,%
\begin{equation*}
\lim_{\varepsilon \rightarrow 0}\mathbb{E}\left[ e^{it\left\langle \tilde{%
\eta}^{\varepsilon },\varphi \right\rangle }\right] =\exp \left[ -\frac{t^{2}%
}{2}\left\langle \partial \varphi ,\partial \varphi \right\rangle \right] ,
\end{equation*}%
where 
\begin{equation}
\left\langle \partial \varphi ,\partial \varphi \right\rangle
=\int_{D}\sum_{\alpha =1}^{d}\left( \frac{\partial \varphi }{\partial
x_{\alpha }}\right) ^{2}dx.  \label{partial}
\end{equation}
\end{theorem}

\begin{remark}
By a standard approximation argument (see, e.g. \cite{M}), one can
strengthen the topology of convergence and replace the space of test
functions in Theorem \ref{dirichletXY} by $H^{\kappa }\left( D\right) ,$ the
Sobolev space of degree $\kappa $, for every $\kappa >4$.
\end{remark}

The argument to prove Theorem \ref{dirichletXY} is fairly robust, and as we
see in the next theorem, it applies to a large class of gradient field
models. We say $\left\{ \theta ^{\varepsilon }\left( \cdot \right) \right\}
\in \mathbb{R}^{D^{\varepsilon }}$ is a zero boundary gradient field on $%
D^{\varepsilon }$ with potential $V:\mathbb{R}\rightarrow \mathbb{R}$, if
its distribution is the Gibbs measure%
\begin{equation}
d\nu _{\beta }=Z_{\beta }^{-1}\exp \left[ -\beta \sum_{\left( i,j\right) \in
\left( D^{\varepsilon }\right) ^{\ast }}V\left( \theta ^{\varepsilon }\left(
i\right) -\theta ^{\varepsilon }\left( j\right) \right) \right] \prod_{i\in
D^{\varepsilon }\backslash \partial D^{\varepsilon }}d\theta ^{\varepsilon
}\left( i\right) \prod_{i\in \partial D^{\varepsilon }}\delta _{0}\left(
d\theta ^{\varepsilon }\left( i\right) \right) ,  \label{grad}
\end{equation}%
where $Z_{\beta }$ is the normalizing constant. For any directed edge $%
b=\left( i,j\right) \in (D^{\varepsilon })^{\ast }$, define the associated
gradient field $\eta ^{\varepsilon }\left( b\right) =\theta ^{\varepsilon
}\left( j\right) -\theta ^{\varepsilon }\left( i\right) $.

\begin{theorem}
\label{gf}Suppose that $\eta ^{\varepsilon }$ is the gradient field
associated with the Gibbs measure~(\ref{grad}), such that $V$ is symmetric, $%
V^{\prime \prime }$ is Lipschitz, and $\inf_{x\in \mathbb{R}}V^{\prime
\prime }\left( x\right) >0$. Also suppose $\beta =\beta \left( \varepsilon
\right) $ satisfies (\ref{temp}). For convenience, we assume that $\tilde{%
\theta}^{\varepsilon }$ (and $\tilde{\eta}^{\varepsilon }$) has been
rescaled by a multiplicative constant so that $V^{\prime \prime }\left(
0\right) =1$. Then the same conclusion as in Theorem \ref{dirichletXY} holds.
\end{theorem}

Theorem \ref{gf} can be proved by the same argument as Theorem \ref%
{dirichletXY}, only replacing the Peierls type estimate (Theorem \ref%
{xypierels}) by Proposition \ref{convexpierels}. Taking $V\left( x\right) =%
\frac{1}{2}x^{2}+\lambda x^{4}$, for $\lambda >0$, Theorem \ref{gf} applies
to anharmonic crystal models (see, e.g., \cite{GK}). Moreover, for gradient
field models with periodic boundary conditions, one can prove a version of
Theorem \ref{gf} for more general potentials $V\left( \cdot \right) $.
Indeed, it suffices to have $V\geq 0$ satisfying the conditions in Theorem %
\ref{gf} in a neighborhood of $0$, and such that $\int_{-\infty }^{\infty
}e^{-bV\left( x\right) }dx<\infty $ for any $b>0$. The proof under this more
general assumption replaces Proposition \ref{convexpierels} by a similar
estimate using reflection positivity \cite{Sh} (see also \cite{Bis}).

Below, when there is no danger of confusion, we will omit some super and
subscripts. For example, we will denote $\tilde{\eta}^{\varepsilon }$ as $%
\eta $, and $\mu _{\beta }^{\varepsilon }$ as $\mu $, etc.

\section{Coupling to convex gradient fields\label{couple}}

In this section, we discuss coupling the XY\ model to convex gradient
fields. The coupling will be used at the end of Section \ref{bdd} to prove
Theorem \ref{dirichletXY}.

\subsection{Contour estimates\label{contour}}

The following estimate for the angular variables in the classical XY model
is due to Bricmont and Fontaine \cite{BF}, the proof of which makes use of
the Ginibre inequality. This is a continuous spin analogue of the Peierls
estimate for discrete spin models.

\begin{theorem}[\protect\cite{BF}]
\label{xypierels}There exists a constant $c<\infty $ such that, for any $%
\varepsilon >0$, any $a\in (0,\pi ]$ and any collection $\mathcal{C}$ of
edges $b\in \left( D^{\varepsilon }\right) ^{\ast }$ such that if $\left(
i,j\right) \in \mathcal{C}$ then $\left( j,i\right) \notin \mathcal{C}$,%
\begin{equation*}
\mathbb{E}\left[ \prod_{b\in \mathcal{C}}1_{\left\vert \eta ^{\varepsilon
}\left( b\right) \right\vert >a}\right] \leq \left( c\exp \left( -\beta
a^{2}/\pi ^{2}\right) \right) ^{\left\vert \mathcal{C}\right\vert }.
\end{equation*}
\end{theorem}

\bigskip Therefore if we fix any $\delta >0$, then as $\beta \rightarrow
\infty $, it is very unlikely to find a large set of edges, such that the
gradient variable $\left\vert \eta ^{\varepsilon }\left( \cdot \right)
\right\vert $ on each of these edges is larger than $\delta $. It is thus
natural to believe that if we modify the cosine potential in (\ref{xygibbs})
outside a $\delta -$neighborhood of $0$, the large scale behavior of the new
gradient field will be the same as for (\ref{xygibbs}), in the $\beta
\rightarrow \infty $ limit.

Given $0<\delta <\pi /2$, define a convex potential $V^{\delta }:\mathbb{%
R\rightarrow R}$ by%
\begin{equation}
V^{\delta }\left( x\right) =\left\{ 
\begin{array}{cc}
-\cos x & \left\vert x\right\vert \leq \delta \\ 
-\cos \delta +\sin \delta \cdot \left( x-\delta \right) +\frac{\cos \delta }{%
2}\left( x-\delta \right) ^{2} & x>\delta \\ 
-\cos \delta -\sin \delta \cdot \left( x+\delta \right) +\frac{\cos \delta }{%
2}\left( x+\delta \right) ^{2} & x<-\delta%
\end{array}%
\right. .  \label{vd}
\end{equation}%
Notice that $\inf_{x\in \mathbb{R}}\left( V^{\delta }\right) ^{\prime \prime
}\left( x\right) =\cos \delta $. Now let $\mu _{\beta }^{\varepsilon ,\delta
}$ be the gradient Gibbs measure defined by 
\begin{equation}
d\mu _{\beta }^{\varepsilon ,\delta }=\left( Z_{\beta }^{\delta }\right)
^{-1}\exp \left[ -\beta \sum_{\left( i,j\right) \in \left( D^{\varepsilon
}\right) ^{\ast }}V^{\delta }\left( \theta ^{\varepsilon ,\delta }\left(
i\right) -\theta ^{\varepsilon ,\delta }\left( j\right) \right) \right]
\prod_{i\in D^{\varepsilon }\backslash \partial D^{\varepsilon }}d\theta
^{\varepsilon ,\delta }\left( i\right) \prod_{i\in \partial D^{\varepsilon
}}\delta _{0}\left( d\theta ^{\varepsilon ,\delta }\left( i\right) \right) .
\label{graddelta}
\end{equation}%
Let $\mathbb{E}^{\delta }$ denote the expectation with respect to $\mu
_{\beta }^{\varepsilon ,\delta }$.

We now state a contour estimate for the measure $\mu _{\beta }^{\varepsilon
,\delta }$. The proof uses the following version of the Brascamp-Lieb
inequality that applies to gradient Gibbs measures with convex potential.

\begin{lemma}[Brascamp-Lieb inequality \protect\cite{BL}]
\label{BL}Let $D$ be a finite graph. Let $\mu $ be a gradient Gibbs measure
on $D$ defined as follows with $V_{e}^{^{\prime \prime }}$ continuous on $%
\mathbb{R}$ for each $e\in D^{\ast }:$ 
\begin{equation}
d\mu =Z^{-1}\exp \left[ -\sum_{e=\left( i,j\right) \in D^{\ast }}V_{e}\left(
\theta \left( i\right) -\theta \left( j\right) \right) \right] \prod_{i\in
D\backslash \partial D}d\theta \left( i\right) \prod_{i\in \partial D}\delta
_{0}\left( d\theta \left( i\right) \right) .
\end{equation}%
For $e=\left( i,j\right) \in D^{\ast }$, let $\eta \left( e\right) =$ $%
\theta \left( i\right) -\theta \left( j\right) $. Let $c_{-}=\inf_{e\in
D^{\ast }}\inf_{x\in \mathbb{R}}V_{e}^{\prime \prime }\left( x\right) $ and
assume $c_{-}>0$. Then the following two inequalities are valid.

1. For all $e\in D^{\ast }$, Var$_{\mu }\left( \eta \left( e\right) \right)
\leq c_{-}^{-1}$.

2. For any $t>0$ and any $\varphi \in \mathbb{R}^{D}$, 
\begin{equation*}
\mathbb{E}\left[ e^{t\left\langle \eta ,\varphi \right\rangle }\right] \leq
\exp \left( \frac{t^{2}}{2}\left\langle \nabla \varphi ,c_{-}^{-1}\nabla
\varphi \right\rangle \right) ,
\end{equation*}%
where $\left\langle \eta ,\varphi \right\rangle =\sum_{e\in D^{\ast }}\eta
\left( e\right) \nabla \varphi \left( e\right) $.
\end{lemma}

\begin{proposition}
\label{convexpierels}For any $\varepsilon >0$, any $0<a<\delta \in (0,\pi
/3] $, any collection $\mathcal{C}$ of edges $b=\left( i,j\right) \in \left(
D^{\varepsilon }\right) ^{\ast }$ such that if $\left( i,j\right) \in 
\mathcal{C}$ then $\left( j,i\right) \notin \mathcal{C}$, and any $\beta
\geq 0$,%
\begin{equation*}
\mathbb{E}^{\delta }\left[ \prod_{\left( i,j\right) \in \mathcal{C}%
}1_{\left\vert \theta ^{\varepsilon ,\delta }\left( j\right) -\theta
^{\varepsilon ,\delta }\left( i\right) \right\vert >a}\right] \leq \left(
\exp \left( 1-\beta a^{2}/\pi ^{2}\right) \right) ^{\left\vert \mathcal{C}%
\right\vert }.
\end{equation*}
\end{proposition}

\begin{proof}
To simplify notation, we write $\theta ^{\varepsilon ,\delta }$ as $\theta $%
, $\eta ^{\varepsilon ,\delta }$ as $\eta $, $\mu _{\beta }^{\varepsilon
,\delta }$ as $\mu $, and $Z^{\delta }$ as $Z$. Note that%
\begin{eqnarray}
&&\mathbb{E}^{\delta }\left[ \prod_{\left( i,j\right) \in \mathcal{C}%
}1_{\left\vert \theta \left( j\right) -\theta \left( i\right) \right\vert >a}%
\right]  \notag \\
&=&\mathbb{E}^{\delta }\left[ \prod_{\left( i,j\right) \in \mathcal{C}%
}1_{\left\vert \theta \left( j\right) -\theta \left( i\right) \right\vert
>a}\exp \left[ \frac{\beta }{2}V^{\delta }\left( \theta \left( j\right)
-\theta \left( i\right) \right) \right] \exp \left[ -\frac{\beta }{2}%
V^{\delta }\left( \theta \left( j\right) -\theta \left( i\right) \right) %
\right] \right]  \notag \\
&\leq &\exp \left[ \frac{\beta }{2}\left( 1-\frac{2a^{2}}{\pi ^{2}}\right)
\left\vert \mathcal{C}\right\vert \right] \mathbb{E}^{\delta }\prod_{\left(
i,j\right) \in \mathcal{C}}\exp \left[ \frac{\beta }{2}V^{\delta }\left(
\theta \left( j\right) -\theta \left( i\right) \right) \right] ,  \label{c1}
\end{eqnarray}%
because $\cos x\leq 1-\frac{2a^{2}}{\pi ^{2}}$ for $0\leq a\leq \left\vert
x\right\vert \leq \pi /3$.

To apply the Brascamp-Lieb inequality (Lemma \ref{BL}) and obtain an upper
bound for 
\begin{equation*}
\mathbb{E}^{\delta }\prod_{\left( i,j\right) \in \mathcal{C}}\exp \left[ 
\frac{\beta }{2}V^{\delta }\left( \theta \left( j\right) -\theta \left(
i\right) \right) \right] ,
\end{equation*}%
we first note that since $V^{\delta }\left( x\right) \leq x^{2}/2-1$ (as can
easily be checked by comparing second derivatives) and hence 
\begin{eqnarray}
\mathbb{E}^{\delta }\prod_{\left( i,j\right) \in \mathcal{C}}\exp \left[ 
\frac{\beta }{2}V^{\delta }\left( \theta \left( j\right) -\theta \left(
i\right) \right) \right] &\leq &\mathbb{E}^{\delta }\prod_{b\in \mathcal{C}%
}\exp \left[ \frac{\beta }{2}\left( \frac{\eta \left( b\right) ^{2}}{2}%
-1\right) \right]  \notag \\
&=&\exp \left( -\frac{\beta }{2}\left\vert \mathcal{C}\right\vert \right) 
\mathbb{E}^{\delta }\prod_{b\in \mathcal{C}}\exp \left[ \frac{\beta }{4}\eta
\left( b\right) ^{2}\right] .  \label{c2}
\end{eqnarray}%
For $\alpha \in \left( 0,\beta /4\right) $, we next define a new gradient
Gibbs measure $\mu _{\alpha }$ (with expectation denoted $\mathbb{E}_{\alpha
}$), such that $d\mu _{\alpha }=\left( Z_{\alpha }\right) ^{-1}\exp \left(
\alpha \sum_{b\in \mathcal{C}}\eta \left( b\right) ^{2}\right) d\mu $, where 
$Z_{\alpha }$ is the corresponding partition function. In other words, for
all $b\in \left( D^{\varepsilon }\right) ^{\ast }$, the new pair potential
is $V_{b}\left( x\right) =\beta V^{\delta }\left( x\right) -\alpha
x^{2}1_{\left\{ b\in \mathcal{C}\right\} }$. Since $\alpha \in \left(
0,\beta /4\right) $, we have%
\begin{equation*}
\inf_{b\in \left( D^{\varepsilon }\right) ^{\ast }}\inf_{x\in \mathbb{R}%
}V_{b}^{^{\prime \prime }}\left( x\right) \geq \beta \cos \delta -\alpha
\geq \beta \cos \frac{\pi }{3}-\frac{\beta }{4}=\beta /4.
\end{equation*}

Then by noting that $\mathbb{E}_{\alpha }\eta \left( b\right) =0$ for any $%
b\in \mathcal{C}$, and applying item 1 of Lemma \ref{BL}, with $c_{-}\geq
\beta /4$ by the last inequality, we obtain%
\begin{eqnarray*}
\log \frac{Z_{\beta /4}}{Z} &=&\int_{0}^{\beta /4}\sum_{b\in \mathcal{C}}%
\mathbb{E}_{\alpha }\left[ \eta \left( b\right) ^{2}\right] d\alpha \\
&=&\int_{0}^{\beta /4}\sum_{b\in \mathcal{C}}\text{Var}_{\alpha }\left[ \eta
\left( b\right) \right] d\alpha \\
&\leq &\int_{0}^{\beta /4}\sum_{b\in \mathcal{C}}\frac{4}{\beta }d\alpha
=\left\vert \mathcal{C}\right\vert ,
\end{eqnarray*}%
or%
\begin{equation}
\mathbb{E}^{\delta }\prod_{b\in \mathcal{C}}\exp \left[ \frac{\beta }{4}\eta
\left( b\right) ^{2}\right] =\frac{Z_{\beta /4}}{Z}\leq e^{\left\vert 
\mathcal{C}\right\vert }.  \label{c3}
\end{equation}%
Combining (\ref{c1})-(\ref{c3}), we have%
\begin{equation*}
\mathbb{E}^{\delta }\left[ \prod_{\left( i,j\right) \in \mathcal{C}%
}1_{\left\vert \theta \left( j\right) -\theta \left( i\right) \right\vert >a}%
\right] \leq \left( e\exp \left( -\beta a^{2}/\pi ^{2}\right) \right)
^{\left\vert \mathcal{C}\right\vert },
\end{equation*}%
which completes the proof.
\end{proof}

In Sections \ref{Zd} and \ref{bdd} we will fix some $\delta \in \left( 0,\pi
/3\right) $, and, abusing the notation when it is clear from the context,
will use $V$, $\mu $, $\mathbb{E}$ to denote the corresponding convex
potential, gradient Gibbs measure and expectation without indicating their
dependence on $\delta $.

\subsection{The Coupling\label{coupling}}

Let $\eta ^{\delta }=\eta ^{\delta ,\varepsilon }$ denote the gradient
variables associated with the Gibbs measure (\ref{graddelta}) on $%
D^{\varepsilon }$ (we will omit its dependence on the parameter $\varepsilon 
$ hereafter). Since $\eta ^{\delta }$ is the gradient of a function, it
satisfies the plaquette condition%
\begin{equation*}
\eta ^{\delta }\left( b_{12}\right) +\eta ^{\delta }\left( b_{23}\right)
+\eta ^{\delta }\left( b_{34}\right) +\eta ^{\delta }\left( b_{41}\right) =0,
\end{equation*}%
for any directed bonds $\left( b_{12},b_{23},b_{34},b_{41}\right) $ that
encloses a plaquette. More generally, $\eta ^{\delta }$ is vortex-free. That
is, for any given directed loop $\mathfrak{L}\subset D^{\varepsilon }$, 
\begin{equation*}
\sum_{b\in \mathfrak{L}}\eta ^{\delta }\left( b\right) =0.
\end{equation*}

Let $\eta $ denote the gradient variables associated with the XY model on $%
D^{\varepsilon }$, defined in (\ref{eta}). Since the spins are $S^{1}-$%
valued, when summing over any directed loop $\mathfrak{L}$, one has 
\begin{equation}
\sum_{b\in \mathfrak{L}}\eta \left( b\right) =2k\pi ,\text{ }k\in \mathbb{Z}%
\text{.}  \label{vortex}
\end{equation}%
When $k\neq 0$, we say $\mathfrak{L}$ is a vortex for $\eta $.

In order to compare the scaling limits under the Gibbs meaures $\mu $
(defined in (\ref{xygibbs})) and $\mu ^{\delta }$ (defined in (\ref%
{graddelta})), we will introduce a coupling between the gradient variables $%
\eta $ and $\eta ^{\delta }$, such that they all agree on an event with
probability close to $1$. In Section \ref{contour} we assume that $\delta
\in (0,\pi /3]$, but for the purposes of the coupling in this section, $%
\delta $ can be anywhere in $\left( 0,\pi /2\right) $. Given some $\omega
\in \mathbb{R}^{\left( D^{\varepsilon }\right) ^{\ast }}$, we will say that
an edge $b\in \left( D^{\varepsilon }\right) ^{\ast }$ is $\delta -$bad (for
that $\omega $), if $\left\vert \omega \left( b\right) \right\vert >\delta $%
; otherwise we say the edge is $\delta -$good. Let $\mathcal{B}$ denote the
subset of $\omega $'s such that there exists at least one $\delta $-bad edge
(for that $\omega $) anywhere in $D^{\varepsilon }$. On $\mathcal{B}^{c}$, $%
\left\vert \eta \left( b\right) \right\vert \leq \delta $ for all $b\in
\left( D^{\varepsilon }\right) ^{\ast }$, and therefore for any directed
bonds $\left( b_{12},b_{23},b_{34},b_{41}\right) $ that encloses a plaquette,%
\begin{equation*}
\left\vert \eta \left( b_{12}\right) +\eta \left( b_{23}\right) +\eta \left(
b_{34}\right) +\eta \left( b_{41}\right) \right\vert \leq 4\delta <2\pi .
\end{equation*}%
Therefore, $\eta \left( b_{12}\right) +\eta \left( b_{23}\right) +\eta
\left( b_{34}\right) +\eta \left( b_{41}\right) =0$. This implies that $\eta 
$ is vortex-free, and for any given loop $\mathfrak{L}$, (\ref{vortex})
cannot hold unless $k=0$, because $D$ is discretely simply connected, as
discussed in Section \ref{model}. Considering the $\eta ^{\delta }$ field as
also taking values in $\mathbb{R}^{\left( D^{\varepsilon }\right) ^{\ast }}$%
, the same considerations apply to it. Thus $\eta $ and $\eta ^{\delta }$
both satisfy the same curl-free constraints on~$\mathcal{B}^{c}$.

Denote by $\left( \partial D^{\varepsilon }\right) ^{\ast }$ the set of
edges in $\left( D^{\varepsilon }\right) ^{\ast }$ with both end vertices in 
$\partial D^{\varepsilon }$. Define 
\begin{equation*}
\Omega =\left\{ \omega :\left( D^{\varepsilon }\right) ^{\ast }\rightarrow 
\mathbb{R}:\sum_{b\in \mathfrak{L}}\omega \left( b\right) =0\text{ for all
loops }\mathfrak{L}\text{, and }\omega \left( b\right) =0\text{ for }b\in
\left( \partial D^{\varepsilon }\right) ^{\ast }\right\} .
\end{equation*}%
First notice that since $V^{\delta }\left( x\right) \geq -\cos x$, we have $%
Z^{\delta }=\mu ^{\delta }\left( \Omega \right) Z^{\delta }\leq \mu \left(
\Omega \right) Z\leq Z$. Let $\mathcal{F}$ be the $\sigma $-field generated
by cylinder sets in $\mathbb{R}^{\left( D^{\varepsilon }\right) ^{\ast }}$.
We now introduce a coupling between $\eta $ and $\eta ^{\delta }$ (thus also
a coupling between $\tilde{\eta}$ and $\tilde{\eta}^{\delta }$) on a common
probability space $\left( \mathbb{R}^{\left( D^{\varepsilon }\right) ^{\ast
}}\times \mathbb{R}^{\left( D^{\varepsilon }\right) ^{\ast }},\mathcal{F}%
\times \mathcal{F}\right) $, with marginals given by $\mu $ and $\mu
^{\delta }$. Let%
\begin{equation*}
\rho =\frac{\mu 1_{\mathcal{B}^{c}}}{\int 1_{\mathcal{B}^{c}}d\mu }=\frac{%
\mu ^{\delta }1_{\mathcal{B}^{c}}}{\int 1_{\mathcal{B}^{c}}d\mu ^{\delta }},
\end{equation*}%
where the second equality follows from the fact that the unnormalized
density of $\mu $ and $\mu ^{\delta }$ coincide on $\mathcal{B}^{c}$ (by the
definition of $V^{\delta }$ in (\ref{vd})). Also let 
\begin{eqnarray*}
\nu &=&\frac{\mu 1_{\mathcal{B}}}{\int 1_{\mathcal{B}}d\mu },\nu ^{\prime }=%
\frac{\mu ^{\delta }1_{\mathcal{B}}}{\int 1_{\mathcal{B}}d\mu ^{\delta }}, \\
c &=&\int 1_{\mathcal{B}^{c}}d\mu ,\text{ }c^{\prime }=\int 1_{\mathcal{B}%
^{c}}d\mu ^{\delta }.
\end{eqnarray*}%
Since on $\mathcal{B}^{c}$ we have $V^{\delta }\left( x\right) =-\cos x$, it
follows that $cZ=c^{\prime }Z^{\delta }$, thus $c\leq c^{\prime }$ since $%
Z^{\delta }\leq Z$. We can therefore write%
\begin{eqnarray*}
\mu &=&c\rho +\left( 1-c\right) \nu \\
\mu ^{\delta } &=&c^{\prime }\rho +\left( 1-c^{\prime }\right) \nu ^{\prime
}=c\rho +\left( 1-c\right) \lambda ,
\end{eqnarray*}%
where $\lambda =\left[ \left( c^{\prime }-c\right) \rho +\left( 1-c^{\prime
}\right) \nu ^{\prime }\right] /\left( 1-c\right) $ is a probability
measure. To describe the coupling, let $\eta ,W$ be independent random
variables taking values in $\mathbb{R}^{\left( D^{\varepsilon }\right)
^{\ast }}$ with distributions $\mu $ and $\lambda $, respectively. Given $%
\eta $, if $\eta \notin \mathcal{B}$, set $\eta ^{\delta }=\eta $;
otherwise, if $\eta \in \mathcal{B}$, define $\eta ^{\delta }=W$. It is
straightforward to check that this $\eta ^{\delta }$ has distribution $\mu
^{\delta }$.

We now check that under this coupling, $\eta =\eta ^{\delta }$ with high
probability. Indeed, 
\begin{equation*}
\mathbb{P}(\eta =\eta ^{\delta })\geq c=\left( 1-\mu \left( \mathcal{B}%
\right) \right) .
\end{equation*}%
Moreover, $\mathcal{B}$ occurs with small probability. Since by Theorem \ref%
{xypierels}, we have%
\begin{equation}
\mu \left( \mathcal{B}\right) \leq \sum_{b\in \left( D^{\varepsilon }\right)
^{\ast }}\mu \left( \left\vert \eta \left( b\right) \right\vert >\delta
\right) \leq 2\varepsilon ^{-d}c\exp \left( -\beta \left( \varepsilon
\right) \delta ^{2}/\pi ^{2}\right) .  \label{bad}
\end{equation}%
For $\delta \in \lbrack \pi /3,\pi /2)$, the right hand side tends to zero
as $\varepsilon \rightarrow 0$ because $\beta \left( \varepsilon \right)
-9d\left\vert \log \varepsilon \right\vert \rightarrow \infty $. Thus under
our coupling, $\mathbb{P}(\eta =\eta ^{\delta })\rightarrow 1$ as $%
\varepsilon \rightarrow 0$.

\begin{remark}
\label{vortexfree}(\ref{bad}) and our earlier discussion show that if $\beta
\left( \varepsilon \right) -9d\left\vert \log \varepsilon \right\vert
\rightarrow \infty $, then with probability close to $1$, for small $%
\varepsilon $ the model is vortex free.
\end{remark}

\section{Fluctuations of gradient fields on $\mathbb{Z}^{d}$\label{Zd}}

In this section we prove the GFF scaling limit for the rescaled gradient
variable asscociated with the gradient Gibbs measure on $\mathbb{Z}^{d}$. As
noted at the end of Section \ref{contour} above, henceforth we use a
simplified notation that does not specify the cutoff parameter $\delta \in
\left( 0,\pi /3\right) $ needed to make $V^{\delta }$ (now denoted $V$)
strictly convex in $\mathbb{R}$. Let $\eta ^{\varepsilon }\in \mathbb{R}%
^{\left( \varepsilon \mathbb{Z}^{d}\right) ^{\ast }}$ be sampled from the
gradient Gibbs measure, defined formally as 
\begin{equation}
d\mu ^{\varepsilon }=Z^{-1}\exp \left[ -\beta \left( \varepsilon \right)
\sum_{b\in \left( \varepsilon \mathbb{Z}^{d}\right) ^{\ast }}V\left( \eta
^{\varepsilon }\left( b\right) \right) \right] \prod_{i\in \varepsilon 
\mathbb{Z}^{d}}d\theta ^{\varepsilon }\left( i\right) ,  \label{zdgrad}
\end{equation}%
such that $\eta ^{\varepsilon }\left( b\right) =\theta ^{\varepsilon }\left(
j\right) -\theta ^{\varepsilon }\left( i\right) $ if $b=\left( i,j\right) $,
and the potential $V:\mathbb{R\rightarrow R}$ satisfies 
\begin{equation}
0<c_{-}\leq V^{\prime \prime }\left( x\right) \leq c_{+}<\infty \text{, }%
V^{\prime \prime }\text{ is Lipschitz continuous, and }V^{\prime \prime
}\left( 0\right) =1.  \label{V}
\end{equation}%
Although it is not a priori clear that the infinite volume Gibbs measure (%
\ref{zdgrad}) is well defined, in \cite{FS} the existence and uniqueness was
proved for the translation invariant and ergodic Gibbs states of (\ref%
{zdgrad}), with any fixed slope (i.e., the macroscopic tilt of the $\theta
^{\varepsilon }$ field). In what follows we focus on the (unique) zero tilt
Gibbs state of (\ref{zdgrad}), which arises as the infinite volume limit of
the zero boundary condition Gibbs measure (\ref{grad}).

Given $\varepsilon >0$, define the rescaled gradient variable $\tilde{\eta}%
^{\varepsilon }=\sqrt{\beta \left( \varepsilon \right) }\eta ^{\varepsilon }$%
. Notice that the law of $\tilde{\eta}^{\varepsilon }$ is given by $\tilde{Z}%
^{-1}\exp \left( -\tilde{H}\left( \tilde{\eta}^{\varepsilon }\right) \right)
\prod_{i\in \varepsilon \mathbb{Z}^{d}}d\tilde{\theta}^{\varepsilon }\left(
i\right) $, with Hamiltonian $\tilde{H}\left( \tilde{\eta}^{\varepsilon
}\right) =\sum_{b\in \left( \varepsilon \mathbb{Z}^{d}\right) ^{\ast }}%
\tilde{V}\left( \tilde{\eta}^{\varepsilon }\left( b\right) \right) $, where 
\begin{equation*}
\tilde{V}\left( x\right) =\tilde{V}_{\varepsilon }\left( x\right) =\beta
\left( \varepsilon \right) V(x/\sqrt{\beta \left( \varepsilon \right) }),
\end{equation*}%
and thus $\tilde{V}^{\prime \prime }\left( x\right) =V^{\prime \prime }(x/%
\sqrt{\beta })$. Moreover,%
\begin{equation}
0<c_{-}\leq \inf_{\varepsilon >0}\inf_{x\in \mathbb{R}}\tilde{V}^{^{\prime
\prime }}\left( x\right) \leq \sup_{\varepsilon >0}\sup_{x\in \mathbb{R}}%
\tilde{V}^{^{\prime \prime }}\left( x\right) \leq c_{+}<\infty .
\label{convex}
\end{equation}%
For test functions $\varphi \in \mathcal{C}_{0}^{\infty }\left( \mathbb{R}%
^{d}\right) $, we define the fluctuation field%
\begin{equation*}
\left\langle \tilde{\eta}^{\varepsilon },\varphi \right\rangle =\varepsilon
^{d/2-1}\sum_{b\in \left( \varepsilon \mathbb{Z}^{d}\right) ^{\ast }}\nabla
\varphi \left( b\right) \tilde{\eta}^{\varepsilon }\left( b\right) .
\end{equation*}%
We now state the following CLT (central limit theorem) type result.

\begin{theorem}
\label{convexGFF}Suppose that the potential $V:\mathbb{R\rightarrow R}$
satisfies (\ref{V}), and $\beta \left( \varepsilon \right) \rightarrow
\infty $ as $\varepsilon \rightarrow 0$. Then for all $t\in \mathbb{R}$, and 
$\varphi \in C_{0}^{\infty }(\mathbb{R}^{d})$,%
\begin{equation*}
\lim_{\varepsilon \rightarrow 0}\mathbb{E}\left[ e^{t\left\langle \tilde{\eta%
}^{\varepsilon },\varphi \right\rangle }\right] =\exp \left[ \frac{t^{2}}{2}%
\left\langle \partial \varphi ,\partial \varphi \right\rangle \right] ,
\end{equation*}%
where $\left\langle \partial \varphi ,\partial \varphi \right\rangle $ is
defined in (\ref{partial}).
\end{theorem}

Notice that when $\beta $ does not depend on $\varepsilon $, the analogue of
Theorem \ref{convexGFF} was proved in \cite{NS}, with the quadratic form of (%
\ref{partial}) replaced by the inverse of an elliptic operator (depending on 
$\tilde{V}$). We will prove Theorem \ref{convexGFF} by adapting the
homogenization argument in \cite{GOS} (which was a probabilistic version of
the argument of \cite{NS}) to the case where $\beta =\beta \left(
\varepsilon \right) $ and $\tilde{V}$ depends on $\varepsilon $.

\subsection{The random walk representation\label{HS}}

The study of scaling limits of the gradient Gibbs measure (\ref{zdgrad}) is
usually based (see \cite{HS}, \cite{DGI}, \cite{GOS}) on a representation of
the covariance of local functions of the field in terms of a random walk in
random conductances, known as the Helffer-Sj{\"{o}}strand representation.

We now introduce this random walk representation, in the simple case where $%
\varepsilon =\beta =1$. It adapts to the general case (where $\beta =\beta
\left( \varepsilon \right) $, and the random walk is on $\varepsilon \mathbb{%
Z}^{d}$) in a straightforward manner. Consider a Langevin dynamics for $%
\theta $, consisting of the system of SDEs, 
\begin{equation}
d\theta _{t}\left( x\right) =-\sum_{b\ni x}V^{^{\prime }}\left( \eta
_{t}\left( b\right) \right) dt+\sqrt{2}dB_{x}\left( t\right) \text{, }x\in 
\mathbb{Z}^{d},  \label{langevin}
\end{equation}%
where $\left\{ B_{x}\left( t\right) \text{, }x\in \mathbb{Z}^{d}\right\} $
is a collection of independent standard Brownian motions. For $b=\left(
x,y\right) \in \left( \mathbb{Z}^{d}\right) ^{\ast }$, set $\eta _{t}\left(
b\right) =\theta _{t}\left( y\right) -\theta _{t}\left( x\right) $. This
Langevin dynamics is reversible with respect to the gradient Gibbs measure (%
\ref{zdgrad}).

The Langevin dyamics induces a Markov dynamics on the gradient field $\eta $%
. We first introduce some notation before describing its generator. For $%
\alpha =1,...,d$, the difference operator $\nabla _{\alpha }$ acting on
functions $f:\mathbb{Z}^{d}\rightarrow \mathbb{R}$ is simply $\nabla
_{\alpha }f\left( x\right) =f\left( x+e_{\alpha }\right) -f\left( x\right) $%
. Its adjoint operator $\nabla _{\alpha }^{\ast }$ is then $(\nabla _{\alpha
}^{\ast }f)\left( x\right) =f\left( x\right) -f\left( x-e_{\alpha }\right) $%
. The generator is given by%
\begin{equation*}
LF\left( \eta \right) =-\sum_{x\in \mathbb{Z}^{d}}\left[ \partial
_{x}^{2}F\left( \eta \right) -\left( \sum_{\alpha =1}^{d}\nabla _{\alpha
}^{\ast }V^{\prime }\left( \eta \left( x,x+e_{\alpha }\right) \right)
\right) \partial _{x}F\left( \eta \right) \right] ,
\end{equation*}%
where $F$ is a smooth local function on $\mathbb{R}^{\left( \mathbb{Z}%
^{d}\right) ^{\ast }}$, and $\partial _{x}=\frac{\partial }{\partial \theta
\left( x\right) }$.

A basic object for probabilistic study of a gradient field model is a random
walk $X\left( t\right) $ on $\mathbb{Z}^{2}$ under the dynamic random
environment $\eta _{t}$. Given the process (\ref{langevin}) with initial
condition $\eta _{0}$, $\left\{ X\left( t\right) \right\} _{t\geq 0}$ is the
random walk that performs nearest neighbor jumps starting from $X\left(
0\right) =x$, with time dependent jump rate $V^{\prime \prime }\left( \eta
_{t}\left( b\right) \right) $ along the directed edge $b$. Since $V$ is
uniformly convex, the random environment is uniformly elliptic. The combined
generator for $\eta _{t}$ and $X\left( t\right) $ is given by $\mathcal{L}%
=L+Q$, where $Q$ is defined by%
\begin{equation*}
QF\left( x,\eta \right) =-\sum_{\alpha =1}^{d}\nabla _{\alpha }^{\ast
}V^{\prime \prime }\left( \eta \left( x,x+e_{\alpha }\right) \right) \nabla
_{\alpha }F\left( x,\eta \right) ,
\end{equation*}%
for any function $F$ that has compact support in $x$, and is a smooth local
function of~$\eta $. Denote by $\mathbb{E}_{x,\eta }$ the quenched
expectation with respect to the law of $X\left( t\right) $, given the
environment process $\left\{ \eta _{t}\right\} _{t\geq 0}$.

We now state the main result for the covariance representation, established
in \cite{GOS}. We set $\partial F\left( x,\eta \right) =\partial _{x}F\left(
\eta \right) $.

\begin{proposition}[\protect\cite{GOS}]
\label{cov}For any twice continuously differentiable and square integrable
local functions $F,G$, such that both $LF$ and $LG$ are $L^{2}$ integrable, 
\begin{equation*}
\mathbb{E}\left[ F\left( \eta \right) G\left( \eta \right) \right] -\mathbb{E%
}\left[ F\left( \eta \right) \right] \mathbb{E}\left[ G\left( \eta \right) %
\right] =\int_{0}^{\infty }\sum_{x\in \mathbb{Z}^{d}}\mathbb{E}\left[
\partial F\left( x,\eta \right) \mathbb{E}_{x,\eta }\left[ \partial G\left(
X\left( t\right) ,\eta \left( t\right) \right) \right] \right] dt,
\end{equation*}%
where $\mathbb{E}$ is the expectation with respect to the Gibbs measure (\ref%
{zdgrad}).
\end{proposition}

\subsection{An averaged CLT}

To apply Proposition \ref{cov}, we replace $\mathbb{Z}^{d},\eta ,V$ by $%
\varepsilon \mathbb{Z}^{d},\tilde{\eta}^{\varepsilon },\tilde{V}$
respectively. The main result in this section is the following.

\begin{proposition}
\label{CLT}As $\varepsilon \rightarrow 0$, $\varepsilon X^{\varepsilon
}\left( t\varepsilon ^{-2}\right) $ converges weakly to the standard
Brownian motion in the Skorohod topology $D\left( \left[ 0,T\right] :\mathbb{%
R}^{d}\right) $, for any $T\geq 0$.
\end{proposition}

The proof of Proposition \ref{CLT} is an adaptation of the argument in
Section 4 of \cite{GOS}, which applies the Kipnis-Varadhan method \cite{KV}.
As we pointed out before, the main difference between the setup in \cite{GOS}
and this paper is that in our setting, the potential $\tilde{V}$ (and thus
the gradient Gibbs measure $\mu ^{\varepsilon }$) changes with $\varepsilon $
because $\beta =\beta \left( \varepsilon \right) $. However, thanks to the
contour estimate (Proposition \ref{convexpierels}) and the stationarity of
the process $\eta _{t}$ in $t$, at any given time $t$, the probability to
have a large $\left\vert \eta _{t}^{\varepsilon }\left( b\right) \right\vert 
$ for any edge $b$ is small. Therefore, at any time $t$ the jump rates of $%
X^{\varepsilon }$ concentrate around $V^{\prime \prime }\left( 0\right) $.
We will use this fact to prove an (averaged) invariance principle for $%
X^{\varepsilon }\left( t\right) $.

Let $\tau _{x}$, $x\in \varepsilon \mathbb{Z}^{d}$ denote the shift by $x$:
for $\eta \in \mathbb{R}^{\left( \varepsilon \mathbb{Z}^{d}\right) ^{\ast }}$%
, $\tau _{x}\eta \left( b\right) =\eta \left( b+x\right) $. For $x=e_{\alpha
}$, $\alpha =1,...,d$, we write $\tau _{\alpha }$ for $\tau _{e_{\alpha }}$.
Also, for any smooth local function $F$ on $\mathbb{R}^{\left( \varepsilon 
\mathbb{Z}^{d}\right) ^{\ast }}$, define $D_{\alpha }F\left( \eta \right)
=F\left( \tau _{\alpha }\eta \right) -F\left( \eta \right) $, and denote by $%
D_{\alpha }^{\ast }$ the adjoint operator of $D_{\alpha }$. Consider the
process that describes the environment seen from the random walk%
\begin{equation*}
\hat{\eta}_{t}^{\varepsilon }=\tau _{-X^{\varepsilon }\left( t\right) }%
\tilde{\eta}_{t}^{\varepsilon }.
\end{equation*}%
Notice that $\hat{\eta}^{\varepsilon }$ is a stationary process in $t$. As
was discussed in \cite{GOS}, its generator is given by a self adjoint
operator%
\begin{equation*}
\mathcal{\tilde{L}}F\left( \eta \right) =LF\left( \eta \right) +\sum_{\alpha
=1}^{d}D_{\alpha }^{\ast }\left[ \tilde{V}^{\prime \prime }\left( \eta
\left( 0,\varepsilon e_{\alpha }\right) \right) D_{\alpha }F\right] \left(
\eta \right) .
\end{equation*}%
We denote by $\mathbb{\tilde{P}}$ and $\mathbb{\tilde{E}}$ the law and
corresponding expectation of the process $\hat{\eta}_{t}^{\varepsilon }$. By
Lemma 4.2 of \cite{GOS}, $\hat{\eta}^{\varepsilon }$ is also time ergodic.

For the proof below, we will switch between two measures: $\mathbb{\tilde{P}}
$ (and $\mathbb{\tilde{E}}$) that characterizes the law of the environment
seen from the random walk $\hat{\eta}_{\cdot }^{\varepsilon }$, and $\mathbb{%
P}$ (and $\mathbb{E}$) that characterizes the stationary measure of $\eta
_{\cdot }^{\varepsilon }$ (i.e., the gradient Gibbs measure (\ref{zdgrad})).

Denote by $\mathcal{F}$ the natural filtration generated by $\left\{ \hat{%
\eta}_{s}^{\varepsilon }\right\} _{s\geq 0}$. Observe that $X^{\varepsilon }$
can be decomposed as an additive functional of the process $\hat{\eta}%
^{\varepsilon }$ plus a martingale with bounded increments: for $\alpha
=1,\dots ,d$, 
\begin{equation}
X_{\alpha }^{\varepsilon }\left( t\right) =M_{\alpha }^{\varepsilon }\left(
t\right) +\int_{0}^{t}\left( D_{\alpha }^{\ast }\tilde{V}^{\prime \prime
}\right) \left( \hat{\eta}_{s}^{\varepsilon }\left( 0,\varepsilon e_{\alpha
}\right) \right) ds,  \label{add}
\end{equation}%
where $M_{\alpha }^{\varepsilon }$ is a $\left( \mathbb{\tilde{P}},\mathcal{F%
}\right) $ martingale, such that%
\begin{equation}
\mathbb{\tilde{E}}\left[ M_{\alpha }^{\varepsilon }\left( t\right) M_{\beta
}^{\varepsilon }\left( t\right) \right] =2t\mathbb{\tilde{E}}\tilde{V}%
^{^{\prime \prime }}\left( \hat{\eta}_{t}^{\varepsilon }\left( 0,\varepsilon
e_{\alpha }\right) \right) \delta _{\alpha \beta }.  \label{q1}
\end{equation}

We now claim $\mathbb{\tilde{E}}\tilde{V}^{^{\prime \prime }}\left( \hat{\eta%
}_{t}^{\varepsilon }\left( 0,\varepsilon e_{\alpha }\right) \right)
\rightarrow 1$ as $\varepsilon \rightarrow 0$. Since $\hat{\eta}%
^{\varepsilon }$ is stationary, and for any $t$, $\eta _{t}^{\varepsilon }$
is distributed according to the translation invariant Gibbs measure (\ref%
{zdgrad}), by Proposition \ref{convexpierels}, there exists $c<\infty $,
such that 
\begin{equation}
\mathbb{P}\left( \left\vert \eta _{t}^{\varepsilon }\left( 0,\varepsilon
e_{\alpha }\right) \right\vert >\beta ^{-1/3}\right) \leq c\exp \left(
-\beta ^{1/3}/\pi ^{2}\right) .  \label{oneedge}
\end{equation}%
On the other hand, when $\left\vert \eta _{t}^{\varepsilon }\left(
0,\varepsilon e_{\alpha }\right) \right\vert \leq \beta ^{-1/3}$, by the
Lipschitz continuity of $V^{\prime \prime }$, there exists $C<\infty $, such
that $\left\vert V^{^{\prime \prime }}\left( \eta _{t}^{\varepsilon }\left(
0,\varepsilon e_{\alpha }\right) \right) -1\right\vert \leq C\beta ^{-1/3}$.
We have, by decomposing the space into the event on the left hand side of (%
\ref{oneedge}) and its complement and using the Lipschitz continuity of $%
V^{\prime \prime }$, 
\begin{eqnarray}
\left\vert \mathbb{\tilde{E}}\tilde{V}^{^{\prime \prime }}\left( \hat{\eta}%
_{t}^{\varepsilon }\left( 0,\varepsilon e_{\alpha }\right) \right)
-1\right\vert &=&\left\vert \mathbb{E}V^{^{\prime \prime }}\left( \eta
_{t}^{\varepsilon }\left( 0,\varepsilon e_{\alpha }\right) \right)
-1\right\vert  \notag \\
&\leq &\left( 1+\sup_{x\in \mathbb{R}}V^{\prime \prime }\left( x\right)
\right) c\exp \left( -\beta ^{1/3}/\pi ^{2}\right) +C\beta ^{-1/3}.
\label{e1}
\end{eqnarray}%
Since $\beta \rightarrow \infty $ as $\varepsilon \rightarrow 0$, we
conclude that $\left\vert \mathbb{\tilde{E}}\tilde{V}^{^{\prime \prime
}}\left( \hat{\eta}_{t}^{\varepsilon }\left( 0,\varepsilon e_{\alpha
}\right) \right) -1\right\vert \rightarrow 0$. Therefore for any $t>0$,%
\begin{equation*}
\mathbb{\tilde{E}}\left[ \left\vert M_{\alpha }^{\varepsilon }\left(
t\right) \right\vert ^{2}\right] \rightarrow 2t\text{ as }\varepsilon
\rightarrow 0\text{.}
\end{equation*}%
By a version of the martingale functional central limit theorem (see for
example \cite{He}), \newline
$\varepsilon M_{\alpha }^{\varepsilon }\left( \varepsilon ^{-2}\cdot \right)
\rightarrow B_{\cdot }$ in the Skorohod topology.

For the second term in (\ref{add}), first notice that if $\beta $ is
independent of $\varepsilon $ (thus $\tilde{V}^{\prime \prime }$ is a fixed
function independent of $\varepsilon $), a central limit theorem for $%
\varepsilon \int_{0}^{t\varepsilon ^{-2}}(D_{\alpha }^{\ast }\tilde{V}%
^{\prime \prime })\left( \hat{\eta}_{s}^{\varepsilon }\left( 0,\varepsilon
e_{\alpha }\right) \right) ds$ follows from Theorem 1.8 and Corollary 1.9 of 
\cite{KV}. Here, because of the dependence of $\tilde{V}^{\prime \prime }$
on $\varepsilon $, we need a triangular array central limit theorem.
However, an analog of Theorem 1.8 of \cite{KV} still holds, stated as Lemma %
\ref{decomp} below. We apply Lemma \ref{decomp} with $\hat{\eta}%
^{\varepsilon }=\hat{\eta}^{\varepsilon }\left( 0,\varepsilon e_{\alpha
}\right) $, and $V^{\varepsilon }=D_{\alpha }^{\ast }\tilde{V}^{\prime
\prime }$. Condition (\ref{l2}) of the lemma is verified by the same
argument as (4.10) of \cite{GOS}, with the constant $C$ only depending on $%
c_{-}$ in (\ref{convex}). It follows that there is for each $\alpha $
another square integrable martingale $\mathcal{M}_{\alpha }^{\varepsilon
}\left( \cdot \right) $, with stationary increments, such that if we define
for $s>0$ and $\alpha =1,\dots ,d$,%
\begin{equation*}
e_{\alpha }^{\varepsilon }\left( s\right) =\int_{0}^{s}\left( D_{\alpha
}^{\ast }\tilde{V}^{\prime \prime }\right) \left( \hat{\eta}%
_{u}^{\varepsilon }\left( 0,\varepsilon e_{\alpha }\right) \right) du-%
\mathcal{M}_{\alpha }^{\varepsilon }\left( s\right) ,
\end{equation*}%
then for any $t>0$, 
\begin{equation}
\lim_{\varepsilon \rightarrow 0}\varepsilon \sup_{s\in \left[ 0,\varepsilon
^{-2}t\right] }\left\vert e_{\alpha }^{\varepsilon }\left( s\right)
\right\vert =0\text{ in probability,}  \label{error}
\end{equation}%
and moreover, the covariance of the increment of $\mathcal{M}_{\alpha
}^{\varepsilon }\left( s\right) $ satisfies%
\begin{eqnarray}
&&\mathbb{\tilde{E}}\left[ \mathcal{M}_{\alpha }^{\varepsilon }\left(
t+1\right) -\mathcal{M}_{\alpha }^{\varepsilon }\left( t\right) \right] %
\left[ \mathcal{M}_{\beta }^{\varepsilon }\left( t+1\right) -\mathcal{M}%
_{\beta }^{\varepsilon }\left( t\right) \right]  \notag \\
&=&2\int_{0}^{\infty }\mathbb{\tilde{E}}\left[ \left( D_{\alpha }^{\ast }%
\tilde{V}^{\prime \prime }\right) \left( \hat{\eta}_{0}^{\varepsilon }\left(
0,\varepsilon e_{\alpha }\right) \right) \left( D_{\beta }^{\ast }\tilde{V}%
^{\prime \prime }\right) \left( \hat{\eta}_{t}^{\varepsilon }\left(
0,\varepsilon e_{\beta }\right) \right) \right] dt.  \label{q2}
\end{eqnarray}

By (\ref{add}) and (\ref{error}), it suffices to prove a martingale central
limit theorem for $\varepsilon (M^{\varepsilon }\left( \varepsilon
^{-2}t\right) +\mathcal{M}_{{}}^{\varepsilon }\left( \varepsilon
^{-2}t\right) )$. To do that we will make use of the fact that this process
has uniformly bounded increments (independent of $\varepsilon $). For fixed $%
\varepsilon >0$, combining (\ref{q1}), (\ref{q2}), and (4.11)--(4.12) of 
\cite{GOS}, we see that the increment of $M^{\varepsilon }+\mathcal{M}%
_{{}}^{\varepsilon }$ from $t$ to $t+1$ has the quadratic form%
\begin{equation*}
q_{\alpha \beta }^{\varepsilon }=2\mathbb{\tilde{E}}\tilde{V}^{^{\prime
\prime }}\left( \hat{\eta}_{0}^{\varepsilon }\left( 0,\varepsilon e_{\alpha
}\right) \right) \delta _{\alpha \beta }-2\int_{0}^{\infty }\mathbb{\tilde{E}%
}\left[ \left( D_{\alpha }^{\ast }\tilde{V}^{\prime \prime }\right) \left( 
\hat{\eta}_{0}^{\varepsilon }\left( 0,\varepsilon e_{\alpha }\right) \right)
\left( D_{\beta }^{\ast }\tilde{V}^{\prime \prime }\right) \left( \hat{\eta}%
_{t}^{\varepsilon }\left( 0,\varepsilon e_{\beta }\right) \right) \right] dt.
\end{equation*}%
This quantity can be characterized through a variational formula on $%
\varepsilon \mathbb{Z}^{d}$ (see, e.g. (2.12) of \cite{GOS}). Thus, arguing
as in \cite{GOS}, we obtain the following analogue of (4.13) of \cite{GOS}:
for any $v\in \mathbb{R}^{d},$%
\begin{equation*}
2\sum_{\alpha =1}^{d}v_{\alpha }^{2}\left[ \mathbb{\tilde{E}}\left( \left[ 
\tilde{V}^{^{\prime \prime }}\left( \hat{\eta}_{0}^{\varepsilon }\left(
0,\varepsilon e_{\alpha }\right) \right) \right] ^{-1}\right) \right]
^{-1}\leq v\cdot q^{\varepsilon }v\leq 2\sum_{\alpha =1}^{d}v_{\alpha }^{2}%
\mathbb{\tilde{E}}\tilde{V}^{^{\prime \prime }}\left( \hat{\eta}%
_{0}^{\varepsilon }\left( 0,\varepsilon e_{\alpha }\right) \right) .
\end{equation*}%
We will prove that as $\varepsilon \rightarrow 0$,%
\begin{equation}
\mathbb{\tilde{E}}\tilde{V}^{^{\prime \prime }}\left( \hat{\eta}%
_{0}^{\varepsilon }\left( 0,\varepsilon e_{\alpha }\right) \right)
\rightarrow 1,\text{ }\mathbb{\tilde{E}}\left( \left[ \tilde{V}^{^{\prime
\prime }}\left( \hat{\eta}_{0}^{\varepsilon }\left( 0,\varepsilon e_{\alpha
}\right) \right) \right] ^{-1}\right) \rightarrow 1,  \label{1}
\end{equation}%
so that 
\begin{equation*}
v\cdot q^{\varepsilon }v\rightarrow 2\sum_{\alpha =1}^{d}v_{\alpha }^{2}.
\end{equation*}%
Therefore the conclusion of the proposition follows from a version of the
functional central limit theorem for martingales with uniformly bounded
increments (see e.g. \cite{He}).

We prove (\ref{1}) by applying Proposition \ref{convexpierels} to show that
as $\varepsilon \rightarrow 0$, the value of $\tilde{V}^{^{\prime \prime
}}\left( \hat{\eta}_{0}^{\varepsilon }\left( 0,\varepsilon e_{\alpha
}\right) \right) $ highly concentrates around $1$. Indeed, $\mathbb{\tilde{E}%
}\tilde{V}^{^{\prime \prime }}\left( \hat{\eta}_{0}^{\varepsilon }\left(
0,\varepsilon e_{\alpha }\right) \right) \rightarrow 1$ follows from (\ref%
{e1}), and a similar argument using (\ref{oneedge}) yields%
\begin{eqnarray*}
\mathbb{\tilde{E}}\left\vert \left[ \tilde{V}^{^{\prime \prime }}\left( \hat{%
\eta}_{0}^{\varepsilon }\left( 0,\varepsilon e_{\alpha }\right) \right) %
\right] ^{-1}-1\right\vert &\leq &\max \left\{ 1-\left( 1+C^{\prime }\beta
^{-1/3}\right) ^{-1},\left( 1-C^{\prime }\beta ^{-1/3}\right) ^{-1}-1\right\}
\\
&&+\left( 1+\left( \inf_{x\in \mathbb{R}}V^{\prime \prime }\left( x\right)
\right) ^{-1}\right) c\exp \left( -\beta ^{1/3}/\pi ^{2}\right) ,
\end{eqnarray*}%
for some $C^{\prime }<\infty $. Since $\beta \left( \varepsilon \right)
\rightarrow \infty $, the right hand side above tends to zero as $%
\varepsilon \rightarrow 0$.

\begin{lemma}
\label{decomp}Let $\hat{\eta}^{\varepsilon }$ be the Markov process with
generator $\mathcal{\tilde{L}}$, defined above following the statement of
Proposition \ref{CLT}. Let $V^{\varepsilon }$ for each $\varepsilon $ be a
function satisfying $\mathbb{\tilde{E}}V^{\varepsilon }\left( \hat{\eta}%
_{0}^{\varepsilon }\right) =0$, and such that there is a constant $C$ (not
depending on $\varepsilon $), so that 
\begin{equation}
\mathbb{\tilde{E}}\left[ V^{\varepsilon }\left( \hat{\eta}_{0}^{\varepsilon
}\right) \phi \left( \hat{\eta}_{0}^{\varepsilon }\right) \right] \leq C%
\mathbb{\tilde{E}}\left[ \phi \left( \hat{\eta}_{0}^{\varepsilon }\right) 
\mathcal{\tilde{L}}\phi \left( \hat{\eta}_{0}^{\varepsilon }\right) \right]
^{1/2}\text{, \ for all }\phi \in \mathcal{D}\left( \mathcal{\tilde{L}}%
\right) .  \label{l2}
\end{equation}%
Define 
\begin{equation*}
Y^{\varepsilon }\left( t\right) =\int_{0}^{t}V^{\varepsilon }\left( \hat{\eta%
}_{s}^{\varepsilon }\right) ds.
\end{equation*}%
Then there exists a square integrable martingale $\mathcal{M}%
_{{}}^{\varepsilon }\left( \cdot \right) $ with stationary increments, such
that for any $t>0$, 
\begin{equation}
\lim_{\varepsilon \rightarrow 0}\sup_{s\in \left[ 0,\varepsilon ^{-2}t\right]
}\varepsilon \left\vert Y^{\varepsilon }\left( s\right) -\mathcal{M}%
_{{}}^{\varepsilon }\left( s\right) \right\vert =0\text{ in probability.}
\label{err1}
\end{equation}%
Moreover, 
\begin{eqnarray}
\mathbb{\tilde{E}}\left\vert \mathcal{M}_{{}}^{\varepsilon }\left(
t+1\right) -\mathcal{M}_{{}}^{\varepsilon }\left( t\right) \right\vert ^{2}
&=&-2\mathbb{\tilde{E}}\left[ V^{\varepsilon }\left( \hat{\eta}%
_{0}^{\varepsilon }\right) \mathcal{\tilde{L}}^{-1}V^{\varepsilon }\left( 
\hat{\eta}_{0}^{\varepsilon }\right) \right]  \notag \\
&=&2\int_{0}^{\infty }\mathbb{\tilde{E}}\left[ V^{\varepsilon }\left( \hat{%
\eta}_{0}^{\varepsilon }\right) V^{\varepsilon }\left( \hat{\eta}%
_{t}^{\varepsilon }\right) \right] dt.  \label{m2m}
\end{eqnarray}

In the case of a discrete time Markov chain we have the following analogous
statement. Let $\left\{ \eta _{j}^{\varepsilon }\right\} $ be a Markov chain
with transition operator $q^{\varepsilon }$, reversible with stationary
distribution $\mu =\mu ^{\varepsilon }$ (and denote by $\mathbb{E}$ the
corresponding expectation). Let $V^{\varepsilon }$ for each $\varepsilon $
be a function satisfying $\mathbb{E}V^{\varepsilon }\left( \eta
_{0}^{\varepsilon }\right) =0$ and such that there is an absolute constant $%
C $, such that 
\begin{equation}
\mathbb{E}\left[ V^{\varepsilon }\left( \eta _{0}^{\varepsilon }\right) \phi
\left( \eta _{0}^{\varepsilon }\right) \right] \leq C\mathbb{E}\left[ \phi
\left( \eta _{0}^{\varepsilon }\right) \left( I-q^{\varepsilon }\right) \phi
\left( \eta _{0}^{\varepsilon }\right) \right] ^{1/2}\text{, \ for all }\phi
\in L^{2}\left( \mu \right) .  \label{l2d}
\end{equation}%
Define%
\begin{equation*}
Y_{n}^{\varepsilon }=\sum_{j=1}^{n}V^{\varepsilon }\left( \eta
_{j}^{\varepsilon }\right) .
\end{equation*}%
Then there exists a square integrable martingale $\left\{ \mathcal{M}%
_{j}^{\varepsilon }\right\} $ with stationary increments, such that for any $%
t>0$,%
\begin{equation}
\lim_{\varepsilon \rightarrow 0}\sup_{1\leq j\leq \left[ \varepsilon ^{-2}t%
\right] }\varepsilon \left\vert Y_{j}^{\varepsilon }-\mathcal{M}%
_{j}^{\varepsilon }\right\vert =0\text{ in probability.}  \label{errd}
\end{equation}%
Moreover,%
\begin{equation}
\mathbb{E}\left\vert \mathcal{M}_{n}^{\varepsilon }-\mathcal{M}%
_{n-1}^{\varepsilon }\right\vert ^{2}=\mathbb{E}\left[ V^{\varepsilon
}\left( \eta _{0}^{\varepsilon }\right) \left( I+q^{\varepsilon }\right)
\left( I-q^{\varepsilon }\right) ^{-1}V^{\varepsilon }\left( \eta
_{0}^{\varepsilon }\right) \right] .  \label{m2md}
\end{equation}
\end{lemma}

\begin{proof}
If we replace $V^{\varepsilon }$ by a function $V$ that does not depend on $%
\varepsilon $, Lemma \ref{decomp} is just a restatement of Theorems 1.3 and
1.8 (for the fixed lattice $\mathbb{Z}^{d}$) of \cite{KV}. Indeed, for fixed 
$\varepsilon >0$, the constructions of the martingales in Theorems 1.3 and
1.8 of \cite{KV} yield (\ref{m2m}) and (\ref{m2md}). Therefore it suffices
to check (\ref{err1}) and (\ref{errd}). As was discussed in \cite{KV}, since
the continuous time processes $Y^{\varepsilon }\left( \cdot \right) $ and $%
\mathcal{M}_{{}}^{\varepsilon }\left( \cdot \right) $ are a.s. right
continous, (\ref{err1}) will follow from the discrete time (\ref{errd})
after replacing the supremum in $\left[ 0,\varepsilon ^{-2}t\right] $ by a
supremum over dyadic times. Therefore we will only sketch the key idea to
prove (\ref{errd}).

The main difference between (\ref{errd}) and Theorem 1.3 of \cite{KV} is
that we need to bound a triangular array $\varepsilon (Y_{j}^{\varepsilon }-%
\mathcal{M}_{j}^{\varepsilon })$. However, the proof of that theorem adapts
to the current setting because (\ref{l2d}) gives control of $V^{\varepsilon
} $ uniformly in $\varepsilon $. As in that theorem, we begin by defining
for all $\varepsilon ,\delta >0,$%
\begin{equation*}
u^{\varepsilon ,\delta }=\left( I-q^{\varepsilon }+\delta \right)
^{-1}V^{\varepsilon },
\end{equation*}%
and two sequences 
\begin{eqnarray*}
\mathcal{M}_{n}^{\varepsilon ,\delta } &=&\sum_{j=0}^{n-1}\left[
u^{\varepsilon ,\delta }\left( \eta _{j+1}^{\varepsilon }\right)
-q^{\varepsilon }u^{\varepsilon ,\delta }\left( \eta _{j}^{\varepsilon
}\right) \right] \\
&=&\sum_{j=0}^{n-1}\left[ u^{\varepsilon ,\delta }\left( \eta
_{j+1}^{\varepsilon }\right) -u^{\varepsilon ,\delta }\left( \eta
_{j}^{\varepsilon }\right) +V^{\varepsilon }\left( \eta _{j}^{\varepsilon
}\right) \right] , \\
e_{n}^{\varepsilon ,\delta } &=&u^{\varepsilon ,\delta }\left( \eta
_{0}^{\varepsilon }\right) -u^{\varepsilon ,\delta }\left( \eta
_{n}^{\varepsilon }\right) .
\end{eqnarray*}%
Then, using the upper bound (\ref{l2d}) for fixed $\varepsilon >0$, the same
argument as in the proof of Theorem 1.3 of \cite{KV} implies the $L^{2}$
limits%
\begin{equation*}
\lim_{\delta \rightarrow 0}\mathcal{M}_{n}^{\varepsilon ,\delta }\dot{=}%
\mathcal{M}_{n}^{\varepsilon },\text{ }\lim_{\delta \rightarrow
0}e_{n}^{\varepsilon ,\delta }\dot{=}e_{n}^{\varepsilon },
\end{equation*}%
such that $Y_{n}^{\varepsilon }=\mathcal{M}_{n}^{\varepsilon
}+e_{n}^{\varepsilon }$. Moreover, the uniform bound (\ref{l2d}) implies the
following analogs of (1.11) -- (1.12) of \cite{KV} along a triagular array:%
\begin{eqnarray}
\lim_{\varepsilon \rightarrow 0}\mathbb{E}\left\vert \mathcal{M}%
_{1}^{\varepsilon ,\varepsilon ^{2}}-\mathcal{M}_{1}^{\varepsilon
}\right\vert ^{2} &=&0,  \label{111} \\
\lim_{\varepsilon _{0},\varepsilon \rightarrow 0}\left\langle u^{\varepsilon
,\varepsilon _{0}}-u^{\varepsilon ,\varepsilon ^{2}},\left( I-q^{\varepsilon
}\right) \left( u^{\varepsilon ,\varepsilon _{0}}-u^{\varepsilon
,\varepsilon ^{2}}\right) \right\rangle &=&0,  \label{111a} \\
\lim_{\varepsilon \rightarrow 0}\varepsilon ^{2}\left\langle u^{\varepsilon
,\varepsilon ^{2}},u^{\varepsilon ,\varepsilon ^{2}}\right\rangle &=&0.
\label{112}
\end{eqnarray}

To prove (\ref{errd}), we write, as in \cite{KV}, 
\begin{equation*}
e_{j}^{\varepsilon }=\left( \mathcal{M}_{j}^{\varepsilon ,\varepsilon ^{2}}-%
\mathcal{M}_{j}^{\varepsilon }\right) +e_{j}^{\varepsilon ,\varepsilon
^{2}}+\varepsilon ^{2}\sum_{k=0}^{j-1}u^{\varepsilon ,\varepsilon
^{2}}\left( \eta _{k}^{\varepsilon }\right) .
\end{equation*}%
It suffices to bound the tail probability for each term. Fix $\alpha >0$. By
Doob's inequality and (\ref{111}),%
\begin{eqnarray*}
\mathbb{P}\left( \sup_{1\leq j\leq \left[ \varepsilon ^{-2}\right]
}\left\vert \mathcal{M}_{j}^{\varepsilon ,\varepsilon ^{2}}-\mathcal{M}%
_{j}^{\varepsilon }\right\vert >\alpha \varepsilon ^{-1}\right) &\leq &\frac{%
\varepsilon ^{2}}{\alpha ^{2}}\mathbb{E}\left\vert \mathcal{M}_{\left[
\varepsilon ^{-2}\right] }^{\varepsilon ,\varepsilon ^{2}}-\mathcal{M}_{%
\left[ \varepsilon ^{-2}\right] }^{\varepsilon }\right\vert ^{2} \\
&\leq &\frac{1}{\alpha ^{2}}\mathbb{E}\left\vert \mathcal{M}%
_{1}^{\varepsilon ,\varepsilon ^{2}}-\mathcal{M}_{1}^{\varepsilon
}\right\vert ^{2} \\
&\rightarrow &0\text{ as }\varepsilon \rightarrow 0.
\end{eqnarray*}%
By (\ref{112}),%
\begin{eqnarray}
\mathbb{P}\left( \sup_{1\leq j\leq \left[ \varepsilon ^{-2}\right]
}\left\vert \varepsilon ^{2}\sum_{k=0}^{j-1}u^{\varepsilon ,\varepsilon
^{2}}\left( \eta _{k}^{\varepsilon }\right) \right\vert >\alpha \varepsilon
^{-1}\right) &\leq &\frac{\varepsilon ^{2}}{\alpha ^{2}}\mathbb{E}\left\vert
\varepsilon ^{2}\sum_{j=0}^{\left[ \varepsilon ^{-2}\right] -1}\left\vert
u^{\varepsilon ,\varepsilon ^{2}}\left( \eta _{j}^{\varepsilon }\right)
\right\vert \right\vert ^{2}  \notag \\
&\leq &\frac{\varepsilon ^{2}}{\alpha ^{2}}\left\langle u^{\varepsilon
,\varepsilon ^{2}},u^{\varepsilon ,\varepsilon ^{2}}\right\rangle  \notag \\
&\rightarrow &0\text{ as }\varepsilon \rightarrow 0.  \label{sumerr}
\end{eqnarray}%
Finally, for fixed $\varepsilon _{0}>0$, we write 
\begin{equation*}
e_{j}^{\varepsilon ,\varepsilon ^{2}}=u^{\varepsilon ,\varepsilon
^{2}}\left( \eta _{0}^{\varepsilon }\right) -u^{\varepsilon ,\varepsilon
^{2}}\left( \eta _{j}^{\varepsilon }\right) =u^{\varepsilon ,\varepsilon
^{2}}\left( \eta _{0}^{\varepsilon }\right) -u^{\varepsilon ,\varepsilon
_{0}}\left( \eta _{j}^{\varepsilon }\right) +\left( u^{\varepsilon
,\varepsilon _{0}}\left( \eta _{j}^{\varepsilon }\right) -u^{\varepsilon
,\varepsilon ^{2}}\left( \eta _{j}^{\varepsilon }\right) \right) .
\end{equation*}%
The supremum over the first two terms on the right hand side above can be
bounded similarly to (\ref{sumerr}). The last term can be bounded by
applying the elementary Lemma 1.4 of \cite{KV} and (\ref{111a}). This yields
(\ref{errd}).
\end{proof}

\subsection{Proof of Theorem \protect\ref{convexGFF}}

Since we have proved Proposition \ref{CLT}, Theorem \ref{convexGFF} will
follow from the approach of \cite{NS} and \cite{BS}, as we sketch below.
First notice that the Hessian matrix $D^{2}\tilde{H}\geq c_{-}I$.

For any $t>0$ and any $\varphi \in \mathcal{C}_{0}^{\infty }\left( \mathbb{R}%
\right) $, set $G^{\varepsilon }\left( t\right) =\mathbb{E}\left[
e^{t\left\langle \tilde{\eta}^{\varepsilon },\varphi \right\rangle }\right] $%
. Applying the exponential Brascamp-Lieb inequality (item 2 of Lemma \ref{BL}%
), we conclude that for $t\in \left[ 0,1\right] $, $\left\{ G^{\varepsilon
}\left( t\right) \right\} _{\varepsilon >0}$ is uniformly bounded. The same
computation as in Section 1.3 of \cite{NS} then yields 
\begin{eqnarray*}
\frac{dG^{\varepsilon }\left( t\right) }{dt} &=&\mathbb{E}\left[
\left\langle \tilde{\eta}^{\varepsilon },\varphi \right\rangle
e^{t\left\langle \tilde{\eta}^{\varepsilon },\varphi \right\rangle }\right]
\\
&=&t\varepsilon ^{d}\mathbb{E}\left[ \left\langle \Delta \varphi ,\mathcal{L}%
^{-1}\Delta \varphi \right\rangle e^{t\left\langle \tilde{\eta}^{\varepsilon
},\varphi \right\rangle }\right] .
\end{eqnarray*}%
By Proposition \ref{2ndmoment} below, we can write 
\begin{equation*}
\frac{dG^{\varepsilon }\left( t\right) }{dt}=t\left\langle \nabla \varphi
,\nabla \varphi \right\rangle G_{\varepsilon }\left( t\right)
+o_{\varepsilon }\left( 1\right) t.
\end{equation*}%
Integrating with respect to $t$ and letting $\varepsilon $ go to zero
completes the proof.

\bigskip

To prove Proposition \ref{2ndmoment}, we need the following $L^{1}-$bound
for the spatial derivatives of heat kernel, established in \cite{DD}.

\begin{lemma}
\label{heatkernel} Let $\eta $ be a stationary elliptic random environment
on $(\mathbb{Z}^{d})^{\ast }$. For $v\in \mathbb{Z}^{d}$, let $X_{\cdot
}^{v} $ be a random walk on $\mathbb{Z}^{d}$ starting at $v$ and with
environment $\eta $. We simply write $X_{\cdot }^{0}$ as $X_{\cdot }$. Then
there exists $C_{1},c_{2}<\infty $ such that for all $t>1$, the discrete
spatial derivative $\nabla _{i}\nabla _{j}$ of the heat kernel $p_{0}^{\eta
}\left( t,x\right) \dot{=}\mathbb{P}\left( X_{t}=x\right) $ can be bounded by%
\begin{equation*}
\mathbb{E}\left\vert \nabla _{i}\nabla _{j}p_{0}^{\eta }\left( t,x\right)
\right\vert \leq \frac{C_{1}}{t^{1+d/2}}e^{-\left\Vert x\right\Vert
^{2}/2c_{2}t}
\end{equation*}%
for $i,j=1,\dots ,d$. Here%
\begin{equation*}
\nabla _{i}\nabla _{j}p_{0}^{\eta }\left( t,x\right) \dot{=}\mathbb{P}\left(
X_{t}^{e_{i}}=x+e_{j}\right) -\mathbb{P}\left( X_{t}^{e_{i}}=x\right) -%
\mathbb{P}\left( X_{t}=x+e_{j}\right) +\mathbb{P}\left( X_{t}=x\right) .
\end{equation*}
\end{lemma}

\begin{proposition}
\label{2ndmoment} For all $\varphi \in \mathcal{C}_{0}^{\infty }(\mathbb{R}%
^{2})$, we have%
\begin{equation*}
\varepsilon ^{d}\left\langle \Delta \varphi ,\mathcal{L}^{-1}\Delta \varphi
\right\rangle \rightarrow \left\langle \nabla \varphi ,\nabla \varphi
\right\rangle \text{ }\ \text{in }L^{2}\text{ as }\varepsilon \rightarrow 0%
\text{.}
\end{equation*}
\end{proposition}

The proof of Proposition \ref{2ndmoment} follows from the same argument as
in Section 4.2 of \cite{BS}, and we omit it here. Roughly speaking, both
sides can be written as a time integral over $[0,\infty )$. Then one can
choose a large $T$, such that the contribution from $[T,\infty )$ is
negligible by the $t^{-\left( 1+d/2\right) }$ decay of the spatial
derivative of the heat kernel (Lemma \ref{heatkernel}). On $\left[ 0,T\right]
$, the desired convergence follows from Proposition \ref{CLT}.

\begin{remark}
If the temperature $\beta ^{-1}$ tends to zero at a much faster rate as $%
\varepsilon \rightarrow 0$, i.e. $\beta \left( \varepsilon \right) \gg
\varepsilon ^{-2d}$, then one can prove a version of Theorem \ref%
{dirichletXY} with periodic boundary conditions, i.e., on a $d$-dimensional
torus $\mathbb{T}_{\varepsilon }^{d}$. This can be done by proving a torus
version of Theorem \ref{convexGFF}, and constructing a coupling similar to
the one in Section \ref{coupling}. Adapting Theorem \ref{convexGFF} to a
torus is straightforward: the dynamic environment is still translation
invariant and stationary, so the homogenization argument still applies.
However, the coupling needs some modification, because on a discrete torus,
there exists macroscopic cycles with nonzero winding number. To ensure $\eta
=\nabla \Phi $, for some function $\Phi :\mathbb{T}_{\varepsilon
}^{d}\rightarrow \mathbb{R}$, $\eta $ cannot have vortices (defined in (\ref%
{vortex})) along these macroscopic cycles. When the temperature is low
enough, i.e. $\beta \left( \varepsilon \right) \gg \varepsilon ^{-2d}$, then
Theorems \ref{xypierels} and \ref{convexpierels} imply that with high
probability, $\left\vert \eta \left( b\right) \right\vert \ll \varepsilon $
for all edges $b\in (\mathbb{T}_{\varepsilon }^{d})^{\ast }$. That rules out
the possibility of vortices along macroscopic cycles. Thus one can couple
the XY model and gradient Gibbs measures in $\mathbb{T}_{\varepsilon }^{d}$
with high probability.
\end{remark}

\section{Bounded domains\label{bdd}}

We now establish a version of Theorem \ref{convexGFF} for gradient fields in
the case of bounded domains with Dirichlet boundary condition. It is also
analogous to Theorem \ref{dirichletXY} which was for XY models in bounded
domains.

\begin{theorem}
\label{dirichletconvex} Suppose that $\eta ^{\varepsilon ,D}$ is the
gradient variable associated with the Gibbs measure~(\ref{grad}). Let $%
\tilde{\eta}^{\varepsilon ,D}=\sqrt{\beta \left( \varepsilon \right) }\eta
^{\varepsilon }$, and for all $\varphi \in C_{0}^{\infty }(D)$, define $%
\left\langle \tilde{\eta}^{\varepsilon ,D},\varphi \right\rangle $ as in (%
\ref{xid}). Assume (\ref{V}) and that $\beta \left( \varepsilon \right)
\rightarrow \infty $ as $\varepsilon \rightarrow 0$. Then 
\begin{equation*}
\lim_{\varepsilon \rightarrow 0}\mathbb{E}\left[ e^{it\left\langle \tilde{%
\eta}^{\varepsilon ,D},\varphi \right\rangle }\right] =\exp \left[ -\frac{%
t^{2}}{2}\left\langle \nabla \varphi ,\nabla \varphi \right\rangle \right] .
\end{equation*}
\end{theorem}

To study the Gibbs measure (\ref{grad}), one can apply a version of the
Helffer-Sj{\"{o}}strand representation described in Section \ref{HS}, and
thus reduce the problem to a random walk in a dynamical random environment,
that is killed on the boundary. However, as opposed to the infinite volume
case, the random environment is no longer translation invariant; thus the
homogenization argument used to prove Proposition \ref{CLT} no longer
applies. This difficulty is resolved in \cite{M} by constructing an
approximate harmonic coupling for the Gibbs measure, as we discuss next. The
approximate harmonic coupling indicates that the Gibbs measure is
"Gaussian-like", because it holds exactly for the discrete Gaussian free
field (DGFF): the law of a DGFF on some domain $D$ with boundary condition $%
f $ is equal in law to that of a zero boundary DGFF on $D$, plus the
discrete harmonic extension of $f$ to~$D$.

We now state the approximate harmonic coupling established in \cite{M}. Let $%
D\subset \varepsilon \mathbb{Z}^{2}$ be a bounded subset that approximates a
smooth simply connected domain. For $x,y\in D$, let dist$\left( x,y\right) $
be the graph distance between $x$ and $y$, and we denote $D\left( r\right)
=\left\{ x\in D:\text{dist}\left( x,\partial D\right) >r\right\} $. The
Ginzburg-Landau measure on $D$ with Dirichlet boundary condition $f$ is
defined by 
\begin{equation}
d\nu ^{f}=Z^{-1}\exp \left[ -\sum_{\left( x,y\right) \in D^{\ast }}V\left(
\theta \left( y\right) -\theta \left( x\right) \right) \right] \prod_{x\in
D\backslash \partial D}d\theta \left( x\right) \prod_{x\in \partial D}\delta
\left( \theta \left( x\right) -f\left( x\right) \right) .  \label{GLD}
\end{equation}

\begin{theorem}[\protect\cite{M}]
\label{1.2}Suppose there exists $\Lambda >0$, such that $f:\partial
D\rightarrow \mathbb{R}$ satisfies $\max_{x\in \partial D}\left\vert f\left(
x\right) \right\vert \leq \Lambda \left\vert \log \varepsilon \right\vert
^{\Lambda }$. Let $\theta $ be sampled from the Ginzburg-Landau measure (\ref%
{GLD}) on $D$ with zero boundary condition, and $\theta ^{f}$ be sampled
from Ginzburg-Landau measure on $D$ with boundary condition $f$. Assume $V$
satisfies (\ref{V}) for some $0<c_{-}<c_{+}<\infty $. Then there exist
constants $c,\gamma ,\delta >0$, that only depends on $c_{-},c_{+}$, so that
if $r>c\varepsilon ^{\gamma }$ then the following holds. There exists a
coupling $\left( \theta ,\theta ^{f}\right) $, such that if $\hat{h}:D\left(
r\right) \rightarrow \mathbb{R}$ is discrete harmonic with $\hat{h}%
|_{\partial D\left( r\right) }=\theta ^{f}-\theta |_{\partial D\left(
r\right) }$, then 
\begin{equation*}
\mathbb{P}\left( \theta ^{f}-\theta \neq \hat{h}\text{ in }D\left( r\right)
\right) \leq c\left( \Lambda \right) \varepsilon ^{\delta }.
\end{equation*}
\end{theorem}

\begin{corollary}
\label{cor}Suppose $V=V_{\varepsilon }$ is such that $\inf_{\varepsilon
>0}c_{-}^{\varepsilon }>0$ and $\sup_{\varepsilon >0}c_{+}^{\varepsilon
}<\infty $, then the same conclusions are valid with $c,\gamma ,\delta $
chosen independent of $\varepsilon $.
\end{corollary}

\begin{proof}
The only difference between Theorem \ref{1.2} and our current setting is
that now $\tilde{V}$ may depend on $\varepsilon $, in a way that 
\begin{equation}
0<c_{-}\leq \inf_{\varepsilon >0}\inf_{x\in \mathbb{R}}\tilde{V}%
_{\varepsilon }^{^{\prime \prime }}\left( x\right) \leq \sup_{\varepsilon
>0}\sup_{x\in \mathbb{R}}\tilde{V}_{\varepsilon }^{^{\prime \prime }}\left(
x\right) \leq c_{+}<\infty .  \label{*}
\end{equation}%
However, by keeping track of the proof in \cite{M}, when $\theta $ and $%
\theta ^{f}$ are sampled from the Gibbs measure (\ref{grad}) with nearest
neighbor potential $\tilde{V}_{\varepsilon }$ that satisfies (\ref{*}), we
claim that one can choose corresponding constants $\left( c,\gamma ,\delta
\right) $ that only depend on $c_{-}$ and $c_{+}$ (and not on $\varepsilon $%
). This claim is implicit in \cite{M}; we sketch here the reasoning. Roughly
speaking, denote by $Q$ the law of $\theta +\hat{h}$, the argument in \cite%
{M} is by controlling the symmetrized relative entropy between $Q$ and $\nu
^{f}$, showing that it can be bounded by $c\left( \Lambda \right)
\varepsilon ^{\delta }$. The symmetrized relative entropy can be written as
a summation over all edges $b\in \left( D^{\varepsilon }\right) ^{\ast }$.
For the edges in the bulk ($b\in \left( D\left( r\right) \right) ^{\ast }$),
one may introduce a coupling between two random walks on random
conductances, arising from the Helffer-Sj{\"{o}}strand representation of the
Gibbs measure with different boundary conditions. The bulk term contributes
an error $c\left( \Lambda \right) \varepsilon ^{\delta _{1}}$, where $\delta
_{1}$ only depends on the Nash continuity constant of the RWRC (thus only
depends on $c_{-},c_{+}$). For the remaining boundary terms, one can gain
some regularity of the field $\theta ^{f}-\theta $ by moving a mesoscopic
distance ($c\varepsilon ^{\gamma }$) away from the boundary, and using an
energy inequality to bound the error by $c\left( \Lambda \right) \varepsilon
^{\delta _{2}}$, where $\gamma $ and $\delta _{2}$ only depend on constants
in the Beurling type estimates for RWRC (thus only depend on $c_{-},c_{+}$).
Therefore, one can take $\delta =\max \left\{ \delta _{1},\delta
_{2}\right\} $, that only depends on $c_{-}$ and $c_{+}$. This completes our
sketch.
\end{proof}

For the rest of this section, we will slightly change the notation, and use
superscripts to distinguish the quantities associated with finite volume
fields. Let $\left\{ \theta ^{\varepsilon ,D}\right\} $ be sampled from the
finite volume Gibbs measure (\ref{grad}), and $\left\{ \eta ^{\varepsilon
}\right\} $ be sampled from the infinite volume gradient Gibbs measure (\ref%
{zdgrad}). Fix an $x_{0}\in \partial D$, and let $x^{\ast }$ be the vertex
in $D^{\varepsilon }$ with minimal distance to $x_{0}$. We construct a Gibbs
measure on $\varepsilon \mathbb{Z}^{d}$ pinned at $x^{\ast }$, denoted as $%
\theta _{0}^{\varepsilon }$, by setting $\theta _{0}^{\varepsilon }\left(
x^{\ast }\right) =0$. For any $x\in \varepsilon \mathbb{Z}^{2}$, take a
chain $\mathcal{C}_{x^{\ast },x}$ connecting $x^{\ast }$ to $x$, and let $%
\theta _{0}^{\varepsilon }\left( x\right) =\sum_{b\in \mathcal{C}_{x^{\ast
},x}}\eta ^{\varepsilon }\left( b\right) $. Define $\tilde{\theta}%
_{0}^{\varepsilon }=\sqrt{\beta \left( \varepsilon \right) }\theta
_{0}^{\varepsilon }$, $\tilde{\theta}^{\varepsilon ,D}=\sqrt{\beta \left(
\varepsilon \right) }\theta ^{\varepsilon ,D}$. We would like to apply
Theorem \ref{1.2} to couple $(\tilde{\theta}_{0}^{\varepsilon },\tilde{\theta%
}^{\varepsilon ,D})$ and thus also couple $\left( \tilde{\eta}^{\varepsilon
},\tilde{\eta}^{\varepsilon ,D}\right) $. To do this, we need to check that
the a-priori bound for $\tilde{\theta}_{0}^{\varepsilon }$ on $\partial
D^{\varepsilon }$ holds with high probability.

\begin{lemma}
\label{7.4}There exists some $\Lambda <\infty $ such that 
\begin{equation*}
\mathbb{P}\left( \max_{x\in \partial D^{\varepsilon }}\left\vert \tilde{%
\theta}_{0}^{\varepsilon }\left( x\right) \right\vert \leq \Lambda
\left\vert \log \varepsilon \right\vert ^{\Lambda }\right) =1-O\left(
\varepsilon ^{8}\right) .
\end{equation*}
\end{lemma}

\begin{proof}
This follows from the same argument as Lemma 7.4 of \cite{M}. That is, since 
$\inf_{\varepsilon >0}\inf_{x}\tilde{V}^{^{\prime \prime }}\left( x\right)
>c_{-}>0$, we can apply the exponential Brascamp-Lieb inequality (Lemma \ref%
{BL}) to control the tail probability of each $\tilde{\theta}%
_{0}^{\varepsilon }\left( x\right) $, and take a union bound. See \cite{M}
for more details.
\end{proof}

We are now in a position to finish the proof of Theorem \ref{dirichletconvex}%
, following the approach in \cite{M}.

\begin{proof}[Proof of Theorem \protect\ref{dirichletconvex}]
Recall that Theorem \ref{convexGFF} implies that for all $\varphi \in
C_{0}^{\infty }\left( D\right) $, $t\in \mathbb{R}$, 
\begin{equation*}
\lim_{\varepsilon \rightarrow 0}\mathbb{E}\left[ e^{it\left\langle \tilde{%
\eta}^{\varepsilon },\varphi \right\rangle }\right] =\exp \left[ -\frac{t^{2}%
}{2}\left\langle \nabla \varphi ,\nabla \varphi \right\rangle \right] .
\end{equation*}

By Corollary \ref{cor} and Lemma \ref{7.4}, on an event $\mathcal{H}%
_{\varepsilon }$ with probability $1-O\left( \varepsilon ^{8}\right)
-O\left( \varepsilon ^{\delta }\right) $, we can apply Corollary \ref{cor}
to construct a coupling of $(\tilde{\theta}^{\varepsilon },\tilde{\theta}%
^{\varepsilon ,D})$. More precisely, there exist $\gamma ,\delta >0$,
independent of $\varepsilon $, such that if we let $\hat{h}^{\varepsilon }$
be the harmonic extension of $\tilde{\theta}_{0}^{\varepsilon }-\tilde{\theta%
}^{\varepsilon ,D}$ from $\partial \left( D^{\varepsilon }\left( \varepsilon
^{\gamma }\right) \right) $ to $D^{\varepsilon }\left( \varepsilon ^{\gamma
}\right) $, on $\mathcal{H}_{\varepsilon }^{c}$ we have $\tilde{\theta}%
^{\varepsilon }-\tilde{\theta}^{\varepsilon ,D}=\hat{h}^{\varepsilon }$ in $%
D^{\varepsilon }\left( \varepsilon ^{\gamma }\right) $, and $\mathbb{P}%
\left( \mathcal{H}_{\varepsilon }^{c}\right) =O\left( \varepsilon ^{\delta
}+\varepsilon ^{8}\right) $.

Given $\varphi \in C_{0}^{\infty }\left( D\right) $, we can choose $%
\varepsilon $ small enough so that supp$\left( \varphi \right) \subset
D^{\varepsilon }\left( \varepsilon ^{\gamma }\right) $. Then on $\mathcal{H}%
_{\varepsilon }$, 
\begin{eqnarray*}
\left\langle \tilde{\eta}^{\varepsilon ,D},\varphi \right\rangle
&=&\left\langle \tilde{\eta}^{\varepsilon },\varphi \right\rangle
-\varepsilon ^{d/2-1}\sum_{b\in \left( D^{\varepsilon }\right) ^{\ast
}}\nabla \varphi \left( b\right) \nabla \left( \tilde{\theta}%
_{0}^{\varepsilon }-\tilde{\theta}^{\varepsilon ,D}\right) \left( b\right) \\
&=&\left\langle \tilde{\eta}^{\varepsilon },\varphi \right\rangle
-\varepsilon ^{d/2-1}\sum_{b\in \left( D^{\varepsilon }\right) ^{\ast
}}\nabla \varphi \left( b\right) \nabla \hat{h}^{\varepsilon }\left( b\right)
\\
&=&\left\langle \tilde{\eta}^{\varepsilon },\varphi \right\rangle ,
\end{eqnarray*}%
where to obtain the last equality, we apply summation by parts and use the
harmonicity of $\hat{h}^{\varepsilon }$. Since $\max \left\{ \left\vert 
\mathbb{E}\left[ e^{it\left\langle \tilde{\eta}^{\varepsilon },\varphi
\right\rangle }1_{\mathcal{H}_{\varepsilon }^{c}}\right] \right\vert
,\left\vert \mathbb{E}\left[ e^{it\left\langle \tilde{\eta}^{\varepsilon
,D},\varphi \right\rangle }1_{\mathcal{H}_{\varepsilon }^{c}}\right]
\right\vert \right\} \leq \mathbb{P}\left( \mathcal{H}_{\varepsilon
}^{c}\right) $, which goes to $0$ as $\varepsilon \rightarrow 0$, this
implies 
\begin{eqnarray*}
\lim_{\varepsilon \rightarrow 0}\mathbb{E}\left[ e^{it\left\langle \tilde{%
\eta}^{\varepsilon ,D},\varphi \right\rangle }\right] &=&\lim_{\varepsilon
\rightarrow 0}\mathbb{E}\left[ e^{it\left\langle \tilde{\eta}^{\varepsilon
,D},\varphi \right\rangle }1_{\mathcal{H}_{\varepsilon }}\right]
=\lim_{\varepsilon \rightarrow 0}\mathbb{E}\left[ e^{it\left\langle \tilde{%
\eta}^{\varepsilon },\varphi \right\rangle }1_{\mathcal{H}_{\varepsilon }}%
\right] \\
&=&\lim_{\varepsilon \rightarrow 0}\mathbb{E}\left[ e^{it\left\langle \tilde{%
\eta}^{\varepsilon },\varphi \right\rangle }\right] \\
&=&\exp \left[ -\frac{t^{2}}{2}\left\langle \nabla \varphi ,\nabla \varphi
\right\rangle \right] .
\end{eqnarray*}
\end{proof}

\begin{proof}[Proof of Theorem \protect\ref{dirichletXY}]
By an abuse of notation, we denote by $\tilde{\eta}_{XY}^{\varepsilon ,D}$
the rescaled gradient variable associated with the XY model in $%
D^{\varepsilon }$. Recall the coupling between the XY\ Gibbs measure (\ref%
{xygibbs}) and the gradient Gibbs measure (\ref{grad}), established in
Section \ref{coupling}. We have%
\begin{equation*}
\left\vert \mathbb{E}\left[ e^{it\left\langle \tilde{\eta}_{XY}^{\varepsilon
,D},\varphi \right\rangle }\right] -\mathbb{E}\left[ e^{it\left\langle 
\tilde{\eta}^{\varepsilon ,D},\varphi \right\rangle }\right] \right\vert
\leq 2\mathbb{P}\left( \tilde{\eta}_{XY}^{\varepsilon ,D}\neq \tilde{\eta}%
^{\varepsilon ,D}\right) \leq 2\mu \left( \mathcal{B}\right) .
\end{equation*}%
Therefore, Theorem \ref{dirichletXY} follows from Theorem \ref%
{dirichletconvex}, and the fact that $\mu \left( \mathcal{B}\right)
\rightarrow 0$ as $\varepsilon \rightarrow 0$.
\end{proof}

\section{Open Questions and Future Work\label{open}}

\subsection{Relaxing the conditions on $\protect\beta \left( \protect%
\varepsilon \right) $}

Theorem \ref{dirichletXY} is in the regime where the temperature $\beta
^{-1} $ goes to zero rapidly as the lattice spacing goes to zero. Note that
by Remark \ref{vortexfree}, $\beta \left( \varepsilon \right) \gg -9d\log
\varepsilon $ ensures that the XY model is "vortex free" in a domain with
diameter $1/\varepsilon $. When $\beta \left( \varepsilon \right)
\rightarrow \infty $ at a slower rate, there will still only be a small
fraction of plaquettes that form vortices. However, as was discussed in \cite%
{KT} and \cite{FSp}, at low temperature, positive vortices (i.e., plaquettes 
$P$ such that $k_{P}>0$ in (\ref{kp})) and negative vortices tend to bind to
each other, like in the dipole gas. This indicates that at low temperature,
macroscopic loops are unlikely to be vortices for the XY model. We therefore
expect the conclusion of Theorem \ref{dirichletXY} to hold as long as $\beta
\left( \varepsilon \right) \rightarrow \infty $ when $\varepsilon
\rightarrow 0$.

\subsection{XY model in an external field}

Given $h>0$, consider the XY model in $D^{\varepsilon }$ with magnetic field 
$h$, given by 
\begin{eqnarray}
d\mu _{\beta ,h}^{\varepsilon } &=&Z_{\beta ,h}^{-1}\exp \left[ \beta
\sum_{\left( i,j\right) \in \left( D^{\varepsilon }\right) ^{\ast }}\cos
\left( \theta ^{\varepsilon }\left( i\right) -\theta ^{\varepsilon }\left(
j\right) \right) +h\sum_{i\in D^{\varepsilon }}\cos \theta ^{\varepsilon
}\left( i\right) \right]  \notag \\
&&\times \prod_{i\in D^{\varepsilon }\backslash \partial D^{\varepsilon
}}d\theta ^{\varepsilon }\left( i\right) \prod_{i\in \partial D^{\varepsilon
}}\delta _{0}\left( d\theta ^{\varepsilon }\left( i\right) \right) .
\label{ext}
\end{eqnarray}

The Gibbs measure (\ref{ext}) fits into the general framework studied by 
\cite{PL}, and their Theorem 2 implies that the spin-spin correlation decays
exponentially with rate at least proportional to $h$. However, the spin wave
picture (see e.g. \cite{MW} and \cite{FSp}) suggests that adding an external
field $h$ is like adding a mass proportional to $\sqrt{h}$ to the GFF, and
the actual decay rate should be proportional to $\sqrt{h}$ for small $h$.
Progress towards this question in $d=3$ has recently been announced \cite%
{BFSnew}. Another step would be to understand the zero temperature limit,
e.g., by proving the following.

\begin{conjecture}
\label{dirichletXY copy(1)}Suppose that $\eta ^{\varepsilon }$ is the
gradient field (as defined in (\ref{eta})) associated with the Gibbs measure
(\ref{ext}), and suppose $\beta =\beta \left( \varepsilon \right) $
satisfies (\ref{temp}). Then for any $t\in \mathbb{R}$ and $\varphi \in
C_{0}^{\infty }\left( D\right) $,%
\begin{equation*}
\lim_{\varepsilon \rightarrow 0}\mathbb{E}\left[ e^{it\left\langle \tilde{%
\eta}^{\varepsilon },\varphi \right\rangle }\right] =\exp \left[ -\frac{t^{2}%
}{2}\left\langle \Delta \varphi ,\left( \Delta _{\sqrt{h}}\right)
^{-1}\Delta \varphi \right\rangle \right] ,
\end{equation*}%
where $\tilde{\eta}^{\varepsilon }=\sqrt{\beta }\eta ^{\varepsilon }$, and
for $m>0$, 
\begin{equation*}
\Delta _{m}=\Delta _{0}+m^{2}.
\end{equation*}%
where $\Delta _{0}$ is the Dirichlet Laplacian on $D$.
\end{conjecture}

When $\beta \left( \varepsilon \right) \gg -9d\log \varepsilon $, one can
still couple the Gibbs measure (\ref{ext}) with a massive gradient Gibbs
measure. The Helffer-Sj{\"{o}}strand representation still applies and leads
to a random walk in a random conductance with exponential killing times. The
main technical difficulty is to establish some kind of massive harmonic
coupling (analogous to Theorem \ref{1.2}) for massive gradient fields.

\bigskip

\paragraph{\textbf{Acknowledgments:}}

We thank Thomas Spencer for useful discussions, S.R.Srinivasa Varadhan for
discussion of \cite{KV}, Aernout van Enter for useful communications and
comments on an earlier version of the paper and Xin Sun for useful
communications. The research of C.M.N. and W.W. was supported in part by
U.S. NSF grants DMS-1007524 and DMS-1507019.

\bibliographystyle{alpha}
\bibliography{GLmax}

\end{document}